\documentclass[11pt]{article}
\usepackage{color}
\usepackage{amsmath}
\usepackage{amssymb}
\usepackage{latexsym}

\usepackage[margin=1in]{geometry}
\usepackage{setspace}
\usepackage{setspace}
\newtheorem{assumption}{Assumption}
\def\qed{ \ \vrule width.2cm height.2cm depth0cm\smallskip}
\newenvironment{proof}{\noindent {\bf Proof.\/}}{$\qed$\vskip 0.1in}

\newcommand{\la}{\langle}
\newcommand{\ra}{\rangle}

\newcommand{\esssup}{\operatornamewithlimits{ess sup}}

\newcommand{\ba}{\begin{array}}
\newcommand{\ea}{\end{array}}
\newcommand{\be}{\begin{equation}}
\newcommand{\ee}{\end{equation}}
\newcommand{\bea}{\begin{eqnarray}}
\newcommand{\eea}{\end{eqnarray}}
\newcommand{\beaa}{\begin{eqnarray*}}
\newcommand{\eeaa}{\end{eqnarray*}}

\def\dbE{\mathbb{E}}
\def\dbF{\mathbb{F}}

\def\dbL{\mathbb{L}}

\def\dbP{\mathbb{P}}
\def\dbR{\mathbb{R}}
\def\dbS{\mathbb{S}}

\def\dbU{\mathbb{U}}
\def\dbV{\mathbb{V}}

\def\a{\alpha}
\def\b{\beta}
\def\g{\gamma}
\def\d{\delta}
\def\e{\varepsilon}

\def\l{\lambda}

\def\si{\sigma}
\def\t{\tau}
\def\f{\varphi}
\def\th{\theta}
\def\o{\omega}

\def\D{\Delta}

\def\L{\Lambda}

\def\O{\Omega}
\def\cA{{\cal A}}
\def\cB{{\cal B}}

\def\cE{{\cal E}}
\def\cF{{\cal F}}

\def\cL{{\cal L}}
\def\cM{{\cal M}}

\def\cP{{\cal P}}

\def\cT{{\cal T}}
\def\cU{{\cal U}}
\def\cV{{\cal V}}

\def\ch{\textsc{h}}
\def\ol{\overline}
\def\ul{\underline}

\def\q{\quad}
\def\qq{\qquad}

\def\pa{\partial}
\def\cd{\cdot}
\def\cds{\cdots}

\def\tr{\hbox{\rm tr}}

\def\qed{ \hfill \vrule width.25cm height.25cm depth0cm\smallskip}

\newcommand{\basa}{\begin{assumption}}
\newcommand{\easa}{\end{assumption}}

\newcommand{\bas}{\begin{assum}}
\newcommand{\eas}{\end{assum}}

\def\limsup{\mathop{\overline{\rm lim}}}

\def\pa{\partial}
\def\wh{\widehat}

 \def\cd{\cdot}
\def\cds{\cdots}

\def\tr{\hbox{\rm tr$\,$}}

\def\dis{\displaystyle}

\def\1{{\bf 1}}

\def\:{\!:\!}
\def\reff#1{{\rm(\ref{#1})}}
\def \proof{{\noindent \bf Proof\quad}}

\def \dbf{{\mathbf{d}}}

at 9pt

\begin{document}

\newtheorem{thm}{Theorem}[section]
\newtheorem{lem}[thm]{Lemma}
\newtheorem{cor}[thm]{Corollary}
\newtheorem{prop}[thm]{Proposition}
\newtheorem{rem}[thm]{Remark}
\newtheorem{eg}[thm]{Example}
\newtheorem{defn}[thm]{Definition}
\newtheorem{assum}[thm]{Assumption}

\renewcommand {\theequation}{\arabic{section}.\arabic{equation}}
\def\thesection{\arabic{section}}

\title{\bf  Pseudo Markovian Viscosity Solutions of Fully Nonlinear Degenerate PPDEs}

\author{Ibrahim  {\sc Ekren}\footnote{ETH Department of Mathematics, ibrahim.ekren@math.ethz.ch.}   
       \and Jianfeng {\sc Zhang}\footnote{University of Southern California, Department of Mathematics, jianfenz@usc.edu. Research supported in part by NSF grant DMS 1413717.}
}\maketitle

\begin{abstract}  In this paper we propose a new type of viscosity solutions for fully nonlinear path dependent PDEs.  By restricting to certain pseudo Markovian structure, we remove the uniform non-degeneracy condition imposed in our earlier works \cite{ETZ1, ETZ2}.  We establish the comparison principle under natural and mild conditions. Moreover, as applications we apply our results to two important classes of PPDEs:   the stochastic HJB equations and  the path dependent Isaacs equations, induced from the stochastic optimization with random coefficients and the path dependent zero sum game problem, respectively.  

\end{abstract}



\vfill\eject

\section{Introduction}
\label{sect-Introduction}
\setcounter{equation}{0}

In this paper we study the  following fully nonlinear parabolic  path-dependent PDE with terminal condition $u(T,\o) = \xi(\o)$: 
 \bea
 \label{PPDE}
 \cL u(t,\o) := \pa_t u(t,\o)
       +G(t,\o,u,\pa_\o u,\pa^2_{\o\o} u)
 =0,~~(t,\o)\in [0, T)\times \O.
 \eea
 Here $\O$ consists of continuous paths $\o$ on $[0, T]$ starting from the origin, $G$ is a progressively measurable generator, and the path derivatives $\pa_t u, \pa_\o u, \pa^2_{\o\o} u$ are defined through a functional It\^{o} formula, initiated by Dupire \cite{Dupire}, see also Cont \& Fournie \cite{CF}. Such equation was first introduced  by Peng {\cite{Peng-ICM, Peng-viscosity}}. In a series of papers Ekren, Keller, Touzi \& Zhang \cite{EKTZ} and Ekren, Touzi \& Zhang   \cite{ETZ1, ETZ2}, we proposed a notion of viscosity solution for such PPDEs and established its wellposedness: existence, comparison principle, and stability. The main innovation of our notion is that, due to the lack of local compactness of the state space $\O$, we replace the pointwise maximum in standard PDE literature with an optimal stopping problem under certain nonlinear expectation.
 
Roughly speaking, the strategy in \cite{ETZ1, ETZ2} is a combination of partial comparison, which is a comparison between a classical semisolution and a viscosity semisolution, and a variation of the Perron's approach. In particular, when the PPDE has a classical solution, it is unique in the sense of viscosity solution, as a direct consequence of the partial comparison.  By utilizing certain path frozen PDE (not PPDE!), in \cite{ETZ2} we established the comparison in the case that the viscosity solution can be approximated by piecewise classical semisolutions taking the form:
\bea
\label{Cbar12}
\sum_{n=0}^\infty v_n\Big((\ch_1,\o_{\ch_1}),\cds, (\ch_n,\o_{\ch_n}); t, \o_t\Big) \1_{\{\ch_n\le t< \ch_{n+1} ~\mbox{or}~\ch_n < \ch_{n+1}=T=t\}},
\eea
where $\ch_n$ is an increasing sequence of stopping times with $\ch_0=0$, and the mapping $(t,x) \mapsto v_n\big((t_1,x_1),\cds, (t_n,x_n); t,x\big)$ is in $C^{1,2}$. However, in order to obtain such smooth $v_n$, we need classical solutions of certain PDEs taking  the form:
\bea
\label{Frozen}
\pa_t v + G_n(t, v, \pa_x v, \pa^2_{xx} v) =0,\q (t,x) \in Q_n\subset [0, T] \times \dbR^d.
\eea
For this purpose in \cite{ETZ2} we have to assume $G$ is uniformly nondegenerate. 

The goal of this paper to is remove this uniform nondegeneracy. We note that degenerate PPDEs appear naturally in many applications, and we will present two examples in this paper. The first one is the stochastic Hamiltonian-Jacobi-Bellman equation, introduced by Peng \cite{Peng-SHJB} to characterize the value function $u(t,x,\o)$ for optimization problems with random coefficients.  \cite{Peng-SHJB} solved the problem when there is only drift control, and the general case with volatility control has been an open problem, see Peng \cite{Peng-open}.  We may view the stochastic HJB equation as a PPDE by considering $x$ as a path. However, this PPDE is by nature always degenerate. We shall characterize the value function as the unique viscosity solution to this degenerate PPDE. We note that in the recent work Qiu \cite{Qiu} viewed the stochastic HJB equation as a backward SPDE and proved its wellposedness in the sense of Sobolev solutions. The second example is the path dependent Isaccs equations, induced from the path dependent zero sum game as  in Pham \& Zhang \cite{PZ2}. In order to obtain the smooth $v_n$ in \reff{Cbar12}, \cite{PZ2} assumes $G$ is uniformly nondegenerate and the dimension $d \le 2$. Besides the degeneracy, our work here also allows for higher dimensions. We note that the recent work Zhang \cite{Zhang} also studied this game problem, and his strategy is in spirit similar to ours. However, no connection with PPDE is studied in \cite{Zhang}.

We still follow the strategy in \cite{ETZ2}, but relying on the viscosity solution theory of PDEs, instead of the classical solution theory of PDEs as in \cite{ETZ2}. Namely, we will construct those $v_n$ via continuous (not $C^{1,2}$~!) viscosity solution of certain path frozen PDE \reff{Frozen}. However, we will establish the uniqueness of viscosity solution in a smaller class. Notice that there is a trade-off between the regularity of the solution and the solution class for the uniqueness. The higher regularity we can establish for the solution (or approximate solution), within a larger solution class we can prove the uniqueness. In our degenerate situation, we are not able to obtain smooth $v_n$ in \reff{Cbar12}, but in $C^0$ only.  As a price, we will establish the uniqueness only in this class, namely there is only one viscosity solution which can be approximated by piecewise Markovian one in the form of  \reff{Cbar12}. This (piecewise) Markovian structure allows us to use the comparison principle of PDE, rather than partial comparison of PPDE. 

There is another major difficulty in the degenerate case. Note that the path frozen PDE \reff{Frozen} is a local PDE, with the domain $Q_n$ induced from the stopping times $\ch_n$. However, in the degenerate case,  the $\ch_n$ used in \cite{ETZ1, ETZ2} has very bad regularity, and consequently the PDE \reff{Frozen} in $Q_n$ typically does not have a continuous viscosity solution. Strongly motivated by the recent work Bayraktar \& Yao \cite{BY}, we shall use some slightly modified stopping times $\ch_n$ which enjoy all the desired properties.   

We remark that the present strategy, as in \cite{ETZ2}, relies heavily on the path frozen PDEs and the related PDE results. In particular, it uses indirectly the very deep regularity results for parabolic PDEs.  In the (possibly degenerate) semi-linear case, Ren, Touzi \& Zhang \cite{RTZ2} and Ren \cite{Ren} studied the regularity for PPDEs directly. The more recent paper Ren, Touzi \& Zhang \cite{RTZ3} established the comparison for fully nonlinear degenerate PPDEs, by introducing certain regularization operator which can be viewed as the counterpart of the sup-convolution in PDE case. Roughly speaking, the strategy in \cite{ETZ2} and the present paper is to approximate the PPDE by certain PDEs and use the solution of the latter to approximate the solution of the original PPDE, while the strategy in \cite{RTZ2, RTZ3} is to approximate the solution of the PPDE directly and show that these approximations are solutions of certain PDEs which are close to the original PPDE in certain sense.  The comparison principle in \cite{RTZ3}, however, is also in a smaller solution class by requiring a somewhat stronger regularity on the solutions, and consequently the coefficients of the PPDE should also have the same stronger regularity.   So there is a tradeoff between \cite{RTZ3} and the present paper: \cite{RTZ3} requires stronger regularity while this paper requires certain piecewise Markovian structure. It will be indeed desirable if one could combine the two techniques and obtain the complete results, which will be left for future research.  

Finally, while we focus on viscosity solutions for PPDEs, there have been different notions of solutions in the literature. First, with the smoothness in terms of Dupire's path derivatives, classical solutions were obtained by Dupire \cite{Dupire} for linear PPDEs (which he called functional PDEs) and by Peng \& Wang \cite{PW} for semilinear PPDEs. Cont \& Fournie \cite{CF} extended the path derivatives to weaker ones  which provided immediately weak solution (in the spirit of Sobolev solution) for linear PPDEs. Peng \& Song \cite{PS} studied Sobolev solution for path dependent HJB equations.  Moreover, Cosso \& Russo \cite{CR} introduced the so called strong-viscosity solution for semilinear PPDEs, which roughly speaking defines the solution as the limit of approximating classical solutions.  While all the notions are consistent with classical solutions when the solutions are smooth,  we emphasize that in path dependent case even the heat equation may not have a classical solution. Our notion of viscosity solution is a local property, and thus the viscosity property  can be easily verified in applications. Of course the challenge lies in the comparison principle, which is the main focus of this paper as well as our earlier works. The Sobolev solution of \cite{PS} is a global solution and involves norm estimates, thus it is easier for uniqueness but more difficult for existence. Indeed, for path dependent Isaacs equations which is a typical example in our approach, it is still not clear what is the appropriate norm under which one may obtain Sobolev solution.  See Pham \& Zhang \cite{PZ1} for some study along this direction.  The strong-viscosity solution of \cite{CR} involves a combination of local and global properties, and is easier for uniqueness but more difficult for existence. Roughly speaking, it transforms the difficulty in our uniqueness to their existence. To our best knowledge, the existence of strong-viscosity solution of \cite{CR} is not clear in fully nonlinear case.

The rest of the paper is organized as follows. In Section 2 we review the basic materials concerning PPDEs. In Section 3 we introduce pseudo Markovian viscosity solutions, and in particular the new hitting times inspired by \cite{BY}. The comparison principle is proved in Section 4, and Section 5 is devoted to existence. In Sections 6 and 7 we present two applications: the stochastic HJB equations induced from  the optimization problem with random coefficients and the path dependent Bellman-Isaacs equations induced from the zero sum stochastic differential games. Finally, some technical proofs are left to Appendix.

\section{Preliminaries}
\label{sect-preliminary}
\setcounter{equation}{0}
 
  In this section, we recall the setup in \cite{ETZ2} and explain why the non-degeneracy requirement is crucial in the uniqueness proof there. 

\subsection{The canonical setting}

Let $\O:= \big\{\o\in C([0,T], \dbR^d): \o_0={\bf 0}\big\}$, the set of continuous paths starting from the origin, $B$ the canonical process, $\dbF = \{\cF_t\}_{0\le t\le T}$ the natural filtration generated by $B$,  $\dbP_0$ the Wiener measure,   $\cT$ the set of $\dbF$-stopping times, and $\L := [0,T]\times \O$. Here and in the sequel, for notational simplicity,  we use ${\bf 0}$ to denote vectors, matrices, or paths with appropriate dimensions whose components are all equal to $0$. 
Moreover, let $\dbS^d$ denote the set of $d\times d$ symmetric matrices, and
 \beaa
 &x \cd x' := \sum_{i=1}^d x_i x'_i 
 ~~\mbox{for any}~~x, x' \in \dbR^d,
 ~~\g : \g' := \tr[\g\g']
 ~~\mbox{for any}~~\g, \g'\in \dbS^d.
 &
 \eeaa
 
 We say a probability measure $\dbP$ on $\cF_T$ is a semimartinagle measure if $B$ is a $\dbP$-semimartingale. For every constant $L>0$, we denote by $\cP_L$ the collection of all semimartingale measures $\dbP$ whose drift and diffusion characteristics are bounded by $L$ and $\sqrt{2L}$, respectively. Denote $\cP_\infty := \cup_{L>0} \cP_L$.

We next discuss regularity of random variables and processes. First, define a semi-norm on $\O$ and a pseudometric on $\L$ as follows: for any $(t, \o), ( t', \o') \in\L$,
\bea\label{rho}
 \|\o\|_{t} 
 := 
 \sup_{0\le s\le t} |\o_s|,
 \q  
 \dbf_\infty\big((t, \o),( t', \o')\big) 
 := 
 |t-t'|^{1\over 2} + \big\|\o_{.\wedge t} - \o'_{.\wedge t'}\big\|_T.
 \eea
For a generic Euclidian space $E$, let  $\dbL^0(\O; E)$ denote the set of   $\cF_T$-measurable random variables $\xi$, $C^0(\O; E)$ (resp. $UC(\O;E)$) the subset of those $\xi$ continuous (resp. uniformly continuous) under $\|\cd\|_T$. Similarly, let  $\dbL^0(\L; E)$ be the set  of $\dbF$-progressively measurable processes $u$, $C^0(\L;E)$  (resp. $UC(\L;E)$) the subset of processes continuous (resp. uniformly continuous) in $(t,\o)$ under $\dbf_\infty$. We use the subscript $_b$  to indicate the subset of bounded elements; and we omit the notation $E$ in the spaces when $E=\dbR$.  For classical solutions of  PPDEs, we need  further regularity of the processes. The following definition  through functional Ito formula is due to \cite{ETZ1} and is inspired by  \cite{Dupire}.
\begin{defn}
\label{defn-spaceC12}  We say $u\in C^{1,2}(\L)$ if $u\in C^0(\L)$ and there exist  $\pa_t u \in C^0(\L)$, $\pa_\o u \in C^0(\L, \dbR^d)$, $\pa^2_{\o\o} u\in C^0(\L, \dbS^d)$ such that, for any $\dbP\in \cP_\infty$, $u$ is a $\dbP$-semimartingale satisfying:
\bea
\label{Ito}
d u = \pa_t u dt+ \pa_\o u \cd d B_t + \frac12 \pa^2_{\o\o} u : d \la B\ra_t,~~0\le t\le T,~~\dbP\mbox{-a.s.}
\eea
\end{defn}
We remark that the  path derivatives  $\pa_t u$, $\pa_\o u$ and $\pa^2_{\o\o} u$, if they exist, are unique.   

We finally introduce the shifted spaces. Let  $0\le s\le t\le T$.

-  Let $\O^t:= \big\{\o\in C([t,T], \dbR^d): \o_t ={\bf 0}\big\}$ be the shifted canonical space and define $B^{t}$, $\dbF^{t}$, $\dbP^t_0$, $\L^t$, $\cT^t$,  $\cP^t_L$, $\cP^t_\infty$ etc. in an obvious sense. In particular, $\L^t :=  [t,T]\times \O^t$.
 Define $\|\cd\|^t_s$ on $\O^t$ and $\dbf^t_\infty$ on $\L^t$ in the spirit of (\ref{rho}),  and the sets $\dbL^0(\L^t; E)$ etc. in an obvious way.   
 
- For  $\o\in \O^s$ and $\o'\in \O^t$, define the concatenation path $\o\otimes_{t} \o'\in \O^s$ by:
\beaa
(\o\otimes_t \o') (r) := \o_r\1_{[s,t)}(r) + (\o_{t} + \o'_r)\1_{[t, T]}(r),
&\mbox{for all}&
r\in [s,T].
\eeaa

- Let $\xi \in \dbL^0(\O^s)$ and $X\in \dbL^0(\L^s)$. For $(t, \o) \in \L^s$, define $\xi^{t,\o} \in \dbL^0(\O^t)$ and $X^{t,\o}\in \dbL^0(\L^t)$ by:
\beaa
\xi^{t, \o}(\o') :=\xi(\o\otimes_t \o'), \q X^{t, \o}(\o') := X(\o\otimes_t \o'),
&\mbox{for all}&
\o'\in\O^t.
\eeaa
It is clear that, for any $(t,\o) \in \L$ and any $u\in C^0(\L)$, we have $u^{t,\o} \in C^0(\L^t)$. The other spaces introduced before enjoy the same property.

\subsection{Viscosity solution of PPDEs}
Our PPDE takes the form of \reff{PPDE} with certain terminal condition $u(T,\o) = \xi(\o)$. We say $u\in C^{1,2}(\L)$ is a classical solution (resp. super-solution, sub-solution) of PPDE \reff{PPDE} if
\beaa
\cL u(t,\o) = ~ (\mbox{resp.}~ \le, \ge )~ 0,\q \forall (t,\o) \in [0, T)\times \O.
\eeaa

The definition of viscosity solution is more involved. First, for any $\xi \in \dbL^0(\O^t)$ with appropriate integrability, we introduce the following  nonlinear expectations:
\bea
\label{cE}
\overline{\cE}^L_t[\xi]
 :=
 \sup_{\dbP\in\cP^t_L}\dbE^{\dbP}[\xi]
 ~~\mbox{and}~~
 \underline{\cE}^L_t[\xi]
 :=
 \inf_{\dbP\in\cP^t_L}\dbE^{\dbP}[\xi]
 =
 -\overline{\cE}^L_t[-\xi].
 \eea
 Next, for any $t\in [0, T]$ and $\e>0$, we define a hitting time:
 \bea
 \label{olch}
 \wh\ch^t_\e := \inf\Big\{s>t: |B^t_s| \ge \e\Big\} \wedge (t+\e) \wedge T.
 \eea
 Now for $u\in \dbL^0(\L)$ with appropriate integrability, we introduce the following classes of test functions: for any $L>0$ and $(t,\o)\in [0, T)\times \O$,
 \bea
 \label{cA}
  \ba{lll}
\dis \underline\cA^{\!L}\!u(t,\o) 
:=
 \Big\{\f\in C^{1,2}(\L^{\!t}):
       (\f-u^{t,\o})_t 
       = 0 =
      \inf_{\t\in \cT^t}\underline\cE^L_t\big[(\f-u^{t,\o})_{\t\wedge \wh\ch^t_\e}
                       \big]
      ~\mbox{for some}~\e>0
 \Big\},
 \\
\dis  \overline\cA^{\!L}\!u(t,\o) 
 :=
 \Big\{\f \in C^{1,2}(\L^{\!t}):
      (\f-u^{t,\o})_t 
      =0=
      \sup_{\t\in \cT^t}\overline \cE^L_t\big[(\f-u^{t,\o})_{\t\wedge\wh\ch^t_\e}
                       \big] 
      ~\mbox{for some}~\e>0
 \Big\}.
 \ea
\eea

\begin{defn}
\label{defn-viscosity}
Let $u\in \dbL^0(\L)$ with appropriate integrability and $L>0$. We say $u$ is a $\cP_L$-viscosity subsolution (resp. supersolution) of PPDE (\ref{PPDE})  if,  for any $(t,\o)\in [0, T)\times \O$ and any $\f \in \underline\cA^{L}u(t,\o)$ (resp. $\f \in \overline\cA^{L}u(t,\o)$):
 \beaa
 \cL^{t,\o}\f(t,{\bf 0}) :=  \pa_t \f(t,{\bf 0})  + G^{t,\o}(\cd, \f,\pa_\o \f,\pa^2_{\o\o}\f) (t, {\bf 0})  &\ge  ~~(\mbox{resp.} \le)&  0.
 \eeaa
We say $u$ is a $\cP_L$-viscosity solution of PPDE (\ref{PPDE}) if it is both a $\cP_L$-viscosity sub- and supersolution.
\end{defn}
We remark that to establish the viscosity theory certain semi-regularity is required for semi-solutions, as introduced in \cite{ETZ1, ETZ2}.  Moreover, the smooth test processes $\f$ in \reff{cA} can actually be restricted to parabolas, and thus the definition can be rewritten in terms of semi-jets, see \cite{RTZ1}.

\subsection{Viscosity solution of PDEs}
In this subsection we consider the following PDE on an open domain $Q \subset [0, T)\times \dbR^d$:
\bea
\label{PDE}
\dbL v(t,x) := \pa_t v (t,x) + g(t,x, v, \pa_x v, \pa_{xx} v) =0,\q (t,x) \in Q.
\eea
We shall introduce two notions of viscosity solutions, one is adapted from Definition \ref{defn-viscosity}, and the other is the standard one in PDE literature, see e.g. \cite{CIL} and \cite{FS}. 

\begin{defn}
\label{defn-viscosityPDE1}
Let $v: Q \to \dbR$ be measurable with certain integrability. 

(i) For some $L>0$, we say $v$ is a $\cP_L$-viscosity subsolution of PDE \reff{PDE}  if,  for any $(t,x)\in Q$,
 \bea
 \label{PDE-cA1}
&\dbL \f(t,x) \ge 0,\q \forall \f\in \ul\cA^L v(t,x), ~\mbox{where}&\\
&\dis \ul\cA^L v(t,x) :=  \big\{\f\in C^{1,2}(Q):  \exists \e>0 ~\mbox{s.t.}~
       [\f-v](t,x) = 0 = \inf_{\t\in \cT^t: \t\le \wh \ch^t_\e }\underline\cE^L_t\big[[\f-v](\t, x+ B^t_{\t}) \big]\big\}.&\nonumber
 \eea
 
 (ii) We say $v$ is a Crandall-Lions viscosity subsolution of PDE \reff{PDE}  if,  for any $(t,x)\in Q$,
  \bea
 \label{PDE-cA2}
 &\dbL \f(t,x) \ge 0,\q \forall \f\in \ul\cA v(t,x), ~\mbox{where}&\\
&\dis \ul\cA v(t,x) :=  \big\{\f\in C^{1,2}(Q): \exists \e>0 ~\mbox{s.t.}~
       [\f-v](t,x) = 0 = \inf_{(s, y)\in Q: |s-t|+|y-x|\le \e}[\f-v](s,y)  \big\}.\nonumber
 \eea

(iii) We define corresponding viscosity supersolution and viscosity solution in an obvious way. 
\end{defn}

\begin{rem}
\label{rem-PDEviscosity}
{\rm (i) When $\e>0$ is small enough, we have $(s, x+ B^t_s) \in Q$ for all $s\le \wh \ch^t_\e$. Thus the $\ul\cA^L v(t,x)$ in \reff{PDE-cA1} is well defined.

(ii) It is clear that $\ul\cA v(t,x)\subset \ul\cA^L v(t,x)$. Then a $\cP_L$-viscosity subsolution has to be a Crandall-Lion viscosity subsolution, but in general not vice versa.

(iii) Because of (ii), formally it could be easier to prove the comparison principle for $\cP_L$-viscosity semi-solutions than for Crandall-Lions viscosity semi-solutions. It will be very interesting to explore such possibility. However, this may require a new argument and we have no clue at this point.
\qed}
\end{rem}

\subsection{The degeneracy of $G$}
\label{sect-degeneracy}
In this subsection we explain why the non-degeneracy requirement is crucial for the comparison principle in \cite{ETZ2}, and how we overcome the difficulties in this paper.

A key element in the strategy of \cite{ETZ2} is the following path frozen PDE (not PPDE!):  for fixed $(t,\o) \in [0, T)\times \O$ and $\e>0$,
\bea
\label{Frozen2}
\pa_t v + G(s, \o_{\cd\wedge t}, v, \pa_x v, \pa_{xx} v) =0,\q (s, x) \in \wh Q^t_\e:= [t, (t+\e)\wedge T) \times \{x\in \dbR^d: |x|< \e\}.
\eea
We emphasize that at above the path $\o$ in $G$ is frozen at $t$ and thus the equation is  a (deterministic) PDE.  Moreover, the domain $\wh Q^t_\e$ is induced by the hitting time $\wh\ch^t_\e$, indeed, we have $(s, B^t_s) \in \wh Q^t_\e$ for $s < \wh\ch^t_\e$.

In order to construct smooth test functions, we let $G_\e$ be a smooth mollifier of $G$ and require the following mollified path frozen PDE (with smooth boundary condition) has a classical solution:
 \bea
\label{Frozen3}
\pa_t v + G_\e(s, \o_{\cd\wedge t}, v, \pa_x v, \pa_{xx} v) =0,\q (s, x) \in \wh Q^t_\e.
\eea
In the PDE literature, one typically needs uniform non-degeneracy of $G_\e$ in terms of $\g$, namely there exists a constant $c_0>0$ such that
\bea
\label{nondegenracy}
G_\e(\cd, \g + \g') - G_\e(\cd, \g) \ge c_0 \tr(\g'),\q \forall \g, \g'\in \dbS^d ~\mbox{with}~\g' \ge {\bf 0}.
\eea
Moreover, for Bellman-Isaacs equations, one may obtain classical solution only when $d\le 2$, even if $G_\e$ is uniformly non-degenerate.   

We note that the classical solution of \reff{Frozen3} is used to prove the partial comparison principle, namely the comparison between a classical semi-solution and a viscosity one.  Our first observation is that, since anyway we are utilizing PDE results, we can use the comparison principle for viscosity solutions of PDE directly. We note that by doing this we are using the regularities of PDEs indirectly, because the comparison principle in PDE literature relies on the regularities through certain regularization procedure. Nevertheless,  this allows us to use the viscosity theory rather that classical theory of PDEs. As a price, however, this requires our viscosity semi-solutions to have certain piecewise Markovian structure, which we will call pseudo Markovian, and thus our comparison principle will be within a smaller class than that in \cite{ETZ2}.

There is another difficulty in degenerate case, even for viscosity theory of PDEs.  Notice that the PDEs \reff{Frozen2} and \reff{Frozen3} are on a bounded domain $\wh Q^t_\e$, not on the whole space. As we see in the following example, in degenerate case such a PDE with smooth boundary condition may not have a continuous viscosity solution.

\begin{eg}
\label{eg-degenerate}
Consider the following degenerate PDE:
\beaa
\pa_t v = 0,~ (t,x) \in  \wh Q^0_\e;\qq v(t,x) = t,~ (t,x) \in \pa \wh Q^0_\e:= \big\{(t,x): t=\e, |x|\le \e~~\mbox{or} ~~ t<\e, |x|=\e\big\}.
\eeaa
Then clearly the candidate solution should be:
$
v(t,x) =\e \1_{\wh Q^0_\e}(t,x) + t \1_{\pa \wh Q^0_\e}(t,x), 
$
which, unfortunately, is discontinuous on $\{(t,x): t=\e, |x|< \e\}$.
\qed
\end{eg}

To overcome this difficulty, we shall modify the hitting time, inspired by \cite{BY}.  While we will study the new hitting time in details in next section, we present a special case here to see how it helps overcome the difficulty in above example. Consider the following hitting time:  
\bea
\label{ch10}
\ch_\e := \inf\{t\ge 0: t + |B_t| \ge \e\},
\eea 
which would induce a domain, changing from a cylinder to a cone:
\bea
\label{O10}
Q_\e := \big\{(t,x) \in [0, \e) \times\dbR^d:  t+ |x|<\e\big\},\q  \pa Q_\e:= \big\{(t,x) \in [0, \e) \times\dbR^d: t + |x| = \e\big\}.
\eea
\begin{eg}
\label{eg-degenerate2}
Consider the following degenerate PDE:
\beaa
\pa_t v = 0,\q (t,x) \in  Q_\e;\qq v(t,x) = t,~ (t,x) \in \pa Q_\e.
\eeaa
Then the  solution is:
$
v(t,x) = \e - |x|, 
$
which is continuous on the whole domain $Q_\e\cup \pa Q_\e$.
\qed
\end{eg}

 \section{Pseudo Markovian viscosity solutions}
 \label{sect-ppde}
\setcounter{equation}{0}

Our PPDE of interest is \reff{PPDE} with terminal condition $u(T,\o) = \xi(\o)$. We shall assume the following standing assumptions.

\begin{assum}
\label{assum-G}
(i)  The PPDE is parabolic, namely $G$ is non-decreasing in $\g$;

(ii) $G$ is uniformly Lipschitz continuous in $(y,z,\g)$ with Lipschitz constant $L$;

(iii) $G$ is  continuous in $(t,\o)$, $G(\cd, 0, {\bf 0}, {\bf 0})$ is bounded, and $\xi\in C^0_b(\O)$. 
\end{assum}
Throughout the paper,  for notational simplicity we denote: for any process $\f$ and $s<t$,
\bea
\label{L1}
L_1:= L+1,\q \f_{s, t} := \f_t - \f_s.
\eea

\subsection{Hitting times}
As explained in Subsection \ref{sect-degeneracy} above, we shall introduce a new type of hitting times, strongly motivated by the recent work  \cite{BY}. Given $R>0$,  $t\in [0, T)$ and $x\in \dbR^d$ with $|x|\le R$, define
\bea
\label{cH}
\ch^{t,x,R}(B^t_{\cdot}) :=   \inf\{s\ge t: |x+B^t_{s}| + L_1(s-t) \ge R\} \wedge T.
\eea
This hitting time enjoys certain nice properties. 
\begin{lem}
\label{lem-cHLip}
For any $(t,x,R)$,  $\t\in \cT^t$ with $\t\le \ch^{t,x,R}$, and $\d>0$, we have
\bea
\label{chDPP}
&\ch^{t,x,R}(B^t_{\cdot}) = \ch^{\t, x+B^t_\t, R-L(\t-t)}(B^t_{\cdot}-B^t_\t),&\\
\label{chlowerbound}
&\sup_{\dbP\in \cP^t_L} \dbP\big(\ch^{t,0,R} < (t + \d)\wedge T\big)  \le C_R\d.&
\eea
Moreover, $\ch^{t,x,R}$ is increasing in $R$, and has the following regularities:
\bea
\label{chLipx}
&\ol\cE^L_t[|\ch^{t,x_1,R_1} - \ch^{t,x_2,R_2}|] \le |x_1-x_2| +|R_1-R_2|,\q |x_1|\le R_1, |x_2|\le R_2;&\\
\label{chLipt}
&\ol \cE^L_{t}[|\ch^{t,x,R}(B^t) - \ch^{\t,x,R}(B^{t}_\cd-B^{t}_{\t})|] \le C\ol\cE^L_t[\sqrt{\t-t}],\q 0\le t \le \t\le \ch^{t,x,R}, |x|\le R.&
\eea
\end{lem}
\proof  First, \reff{chDPP} and the monotonicity of $\ch^{t,x,R}$ in $R$ are obvious. Next, for any $\d>0$, if $L_1 \d \ge {R\over 2}$, then \reff{chlowerbound} becomes trivial. Now assume $L_1 \d \le {R\over 2}$. For any $\dbP\in \cP^t_L$, 
\beaa
 \dbP\big(\ch^{t,0,R} < (t + \d)\wedge T\big)  &\le& \dbP\big(\sup_{t\le s\le t+\d} |B^{t}_s|  + L_1\d  \ge R\big) \le \dbP\big(\sup_{t\le s\le t+\d} |B^{t}_s|   \ge {R\over 2} \big) \\
&\le&{4\over R^2} \dbE^\dbP\Big[\sup_{t\le s\le t+\d} |B^{t}_s|^2\Big] \le C_R \d.
\eeaa
By the arbitrariness of $\dbP$, this implies \reff{chlowerbound}. 
Moreover,  \reff{chLipt} follows directly from \reff{chDPP}, \reff{chLipx},  and the following simple estimate: 
$
\ol \cE^L_t[|B^t_\t|] \le C\ol \cE^L_t[\sqrt{\t-t}].
$

To prove  \reff{chLipx}, we assume without loss of generality that $t=0$ and denote  $\t_i:= \ch^{0,x_i, R_i}$, $i=1,2$, and $\D \f := \f_2-\f_1$ for $\f = x, R, \t$. On $\{\t_1<\t_2\}\in \cF_{\t_1}$ and under each $\dbP\in \cP_L$,  we have
\beaa
&& |x_1 + B_{\t_1}|+L_1 \t_1 = R_1,  \q  |x_2 + B_{\t_2}| +L_1 \t_2 \le R_2\\
&\Longrightarrow& \D R \ge  \dbE^\dbP_{\t_1}\big[|x_2 + B_{\t_2}| + L_1 \t_2\big] - \big[|x_1 + B_{\t_1}|+L_1 \t_1\big]\\
&&\qq\ge \big|x_2 + \dbE^\dbP_{\t_1}[B_{\t_2}]\big|  -  |x_1 + B_{\t_1}| + L_1 \dbE^\dbP_{\t_1}[\D\t] \ge  L_1 \dbE^\dbP_{\t_1}[\D\t] - |\D x| -  \big|\dbE^\dbP_{\t_1}[B_{\t_1, \t_2}]\big|\\
&&\qq \ge L_1 \dbE^\dbP_{\t_1}[\D\t] - |\D x| -  L\dbE^\dbP_{\t_1}[\D\t]  =  \dbE^\dbP_{\t_1}[\D\t] - |\D x| \\
&\Longrightarrow& \dbE^\dbP_{\t_1}[\D\t]\le |\D x| + |\D R|.
\eeaa
This implies that
  \beaa
  \dbE^\dbP\big[(\t_2-\t_1)\1_{\{\t_1 < \t_2\}}\big] \le \big[|\D x| + |\D R|\big]\dbP(\t_1 < \t_2).
  \eeaa
  Similarly, we have
$  \dbE^\dbP\big[(\t_1-\t_2)\1_{\{\t_2 < \t_1\}}\big] \le   \big[|\D x| + |\D R|\big] \dbP(\t_2 < \t_1).$
  Then
$
   \dbE^\dbP[|\D\t|] \le   \big[|\D x| + |\D R|\big],
  $
and \reff{chLipx} follows from the arbitrariness of $\dbP\in \cP_L$.
\qed

\begin{rem}
\label{rem-hitting}
{\rm (i) The regularities \reff{chLipx} and \reff{chLipt} are in $\dbL^1$-sense. The hitting time $\wh\ch^t_\e$ in \reff{olch} shares these properties in uniformly non-degenerate case, but does not  in degenerate case.  The work \cite{BY} introduced a different hitting time which has stronger regulairty: 
\bea
\label{BYch}
\ch^*_\e := \inf\big\{t\ge 0: t + \sup_{0\le s\le t} |B_s| \ge \e\big\}.
\eea
One can easily show that $\ch^*_\e$ is Lipschitz continuous in $\o$ in pathwise sense: 
\beaa
| \ch^*_\e (\o) - \ch^*_\e(\tilde \o)|\le \|\o - \tilde\o\|_T.
\eeaa
However,  $\ch^*_\e$ does not share the Markovian property in the sense of \reff{chDPP}:
\beaa
\ch^*_\e \neq \ch^{*,\t, B_\t}_\e ~\mbox{for}~\t < \wh\ch_\e,\q\mbox{where}~ \wh \ch^{*,t,x}_\e := \inf\big\{s\ge t: s + \sup_{t\le r\le s} |x+B^t_r| \ge \e\big\}.
\eeaa
In this paper we need both the regularity and the Markovian structure, in order to utilize the viscosity theory of PDEs. 

(ii) The regularities \reff{chLipx} and \reff{chLipt}  are under nonlinear expectation. Under standard (linear) expectation, such regularities have been well understood, see e.g.   \cite{Mikulevicious, MR}.

(iii) For any $\e>0$, there exist $0<\e_1, \e_2 < \e$ such that 
\bea
\label{chequivalent}
\wh \ch^t_{\e_1} \le \ch^{t,{\bf 0}, \e},\q  \ch^{t,{\bf 0}, \e_2} \le \wh \ch^t_{\e}.
\eea
Then clearly Definition \ref{defn-viscosity} remains equivalent if we replace the $\wh\ch^t_\e$ in \reff{cA} with $\ch^{t,{\bf 0}, \e}$. Moreover, the optimal stopping problem, which is required in \cite{ETZ1, ETZ2} and proved in \cite{ETZ0}, becomes a lot easier if we use $\ch^{t,{\bf 0}, \e}$ due to the regularities in Lemma \ref{lem-cHLip}. 
\qed}
\end{rem}

The next property will be crucial to pass the local structure to a global one. Fix $\e>0$, define 
\bea
\label{che}
\ch^\e_0 := 0,\q \ch^\e_{n+1} := {\ch^{\ch_n^\e, 0, \e}(B_\cd-B_{\ch^\e_n})} = \inf\{t\ge \ch^\e_n: |B_{\ch^\e_n,t}| + L_1(t - \ch^\e_n) \ge \e\} \wedge T,~ n\ge 0.
\eea
\begin{lem}
\label{lem-chn}
For any $\e>0$,
\bea
\label{Echn}
\bigcap_{n\ge 1} \{\ch^\e_n <T\} = \emptyset\q\mbox{and}\q   \sup_{\dbP\in \cP_L} \dbP [\ch^\e_n< T] \le {C\over n\e^2}.
\eea
 \end{lem}
 \proof First, clearly $\ch^\e_n$ is nondecreasing, and thus $\ch^\e_\infty:= \lim_{n\to\infty} \ch^\e_n\le T$ exists. Note that,
 \beaa
 \bigcap_{n\ge 1} \{\ch^\e_n <T\} \subset \bigcap_{n\ge 1}  \Big\{|B_{\ch^\e_n,\ch^\e_{n+1}}| + L_1(\ch^\e_{n+1} - \ch^\e_n) = \e\Big\}.
 \eeaa
 Since $\lim_{n\to\infty} B_{\ch^\e_n} = B_{\ch^\e_\infty}$, clearly the right side above is empty, then so is the left side.
 
 Next, for any  $n\ge 1$,
 \beaa
 \{\ch^\e_n <T\} \subset \bigcap_{0\le k< n}\Big\{|B_{\ch^\e_k,\ch^\e_{k+1}}|+L_1(\ch^\e_{k+1}- \ch^\e_k) = \e\Big\}
 \eeaa 
 Note that
 \beaa
 \Big(|B_{\ch^\e_k,\ch^\e_{k+1}}|+L_1(\ch^\e_{k+1}- \ch^\e_k)\Big)^2 \le  2|B_{\ch^\e_k,\ch^\e_{k+1}}|^2 + C(\ch^\e_{k+1}- \ch^\e_k).
 \eeaa
 Then
 \beaa
  \{\ch^\e_n <T\} &\subset& \bigcap_{0\le k< n}\Big\{2|B_{\ch^\e_k,\ch^\e_{k+1}}|^2 + C(\ch^\e_{k+1}- \ch^\e_k) \ge \e^2\Big\}\\
   & \subset& \Big\{\sum_{k=0}^{n-1}\big[2|B_{\ch^\e_k,\ch^\e_{k+1}}|^2 + C(\ch^\e_{k+1}- \ch^\e_k)\big] \ge  n\e^2\Big\} \subset \Big\{2\sum_{k=0}^{n-1}|B_{\ch^\e_k,\ch^\e_{k+1}}|^2 + C\ge  n\e^2\Big\}.
 \eeaa 
 Now for  any $\dbP\in\cP_L$,
\beaa
\dbP(\ch^\e_n <T) \le {1\over n\e^2}\dbE^\dbP\Big[2\sum_{k=0}^{n-1}|B_{\ch^\e_k,\ch^\e_{k+1}}|^2 + C\Big]  \le {C\over n\e^2}.
\eeaa
By the arbitrariness of $\dbP\in\cP_L$ we obtain  the second claim immediately.
\qed

\subsection{Piecewise Markovian processes}

For any $t\in [0, T]$ and $\e>0$, denote
 \bea
 \label{Oet}
 \left.\ba{c}
 \dis Q^{\e}_t :=  \big\{(s,x)\in (t, T]\times \dbR^d: |x| + L_1(s-t) <\e\big\};\\
{\dis  \pa Q^{\e}_t 
 :=\big\{(s,x)\in (t, T]\times \dbR^d: |x| + L_1(s-t) =\e\big\} ~\cup~ \big\{(T,x):  |x| + L_1(T-t) \le \e\big\} ; }\\
 {\dis \wh Q^\e_t := Q^\e_{t} \cup \pa Q^\e_{t} \cup \{(t, {\bf 0})\};}\\
\dis  \Pi^\e_n := \big\{\pi_n = (t_i, x_i)_{1\le i\le n}: 0=t_0<t_1<\cds<t_n<T,   (t_i, x_i) \in \pa Q^\e_{t_{i-1}}, 1\le i\le n\big\};\\
\dis D^\e_{n+1} := \big\{(\pi_n; t,x):  \pi_n \in \Pi^\e_n, (t,x) \in \wh Q^\e_{t_n} \big\};\\
\dis \pi^\e_n(\o) := (\ch^\e_i(\o), \o_{\ch^\e_i(\o)} -\o_{\ch^\e_{i-1}(\o)} )_{1\le i\le n}.
\ea\right.
 \eea

In light of \reff{Cbar12}, we introduce processes with the following piecewise Markovian structure.

\begin{defn}
\label{defn-Markov}
Let $\e>0$. We say a process $u\in \dbL^0(\L)$ is $\e$-Markovian, denoted as $u\in \cM_\e(\L)$, if there exist  deterministic functions $v_n: D^\e_{n+1} \to \dbR$, $n\ge 0$, satisfying: 

(i) \reff{Cbar12} holds, namely
\bea
\label{Markov}
u (t,\o)=  \sum_{n=0}^\infty v_n\left(\pi^\e_n(\o); t, \o_t-\o_{\ch_n (\omega)} \right) \1_{\{\ch^\e_n\le t< \ch^\e_{n+1} ~\mbox{or}~\ch^\e_n < \ch^\e_{n+1}=T=t\}}
\eea

(ii) For all $\pi_n = (t_i, x_i)_{1\le i\le n}\in \Pi^\e_n$ and $(t,x)\in \pa Q^\e_{t_n}$,  the following compatibility condition holds 
\bea
\label{compatibility}
v_n(\pi_n; t, x)= v_{n+1}(\pi_n, (t,x); t,0).
\eea

(iii) Each $v_n$, $n\ge 0$, is continuous in $D^\e_{n+1}$.
\end{defn}

\begin{rem}
\label{rem-Markov}
{\rm (i) The continuity of $v_n$ and compatibility  \reff{compatibility}  imply that $u$ is continuous in time. 

(ii) We do not require $v_n(\pi_n; \cd)$ to be continuous on  $\{(t_n, x): 0<|x|\le \e\}$. 
 However, for any $\d>0$ small, $v_n$ is  (uniformly) continuous on the following compact set 
\bea
\label{continuity}
D^{\e,\d}_{n+1}:=\Big\{(\pi_n; t,x) \in D^\e_{n+1}: t_i-t_{i-1}\ge \d, i=1,\cds, n,~\mbox{and}~ t-t_n\ge \d\Big\}.
\eea
It turns out that this  uniform continuity and the continuity of $v_n(\pi_n;\cdot)$ at $(t_n,{\bf 0})$ is enough for our comparison result.  

(iii) In \cite{ETZ2} we imposed a technical condition Assumption 3.5 to ensure the constructed $v_n$ will be uniformly continuous in $D^\e_{n+1}$. This condition is not needed here because of the introduction of our new hitting time. 
\qed}\end{rem}

Moreover, we may extend all the notations to the shifted spaces: given $0\le t<T$,
\bea
\label{Mt}
\ch^{t,\e}_n,\q \Pi^{t,\e}_n,\q \pi^{t,\e}_n,\q \cM_\e(\L^t),~ \mbox{etc.}
\eea

\subsection{Pseudo Markovian viscosity solution}

We provide the following notion of viscosity solutions.

\begin{defn}
\label{defn-viscosity2}
We say $u$ is a pseudo Markovian $\cP_L$-viscosity sub-solution (resp. Crandall-Lions viscosity sub-solution) of PPDE \reff{PPDE}  at $(t,\o) \in [0, T)\times \O$ if, for any $\e>0$, there exists $u^{t, \o, \e} \in \cM_\e(\L^t)$ with corresponding $\{v_n, n\ge 1\}$, such that

(i)  for each $\pi_n = (t_i, x_i)_{1\le i\le n} \in \Pi^{t,\e}_n$, $v_n(\pi_n; \cd)$ is a $\cP_L$-viscosity sub-solution (resp. Crandall-Lions viscosity sub-solution)  to the following PDE: 
\bea
\label{vnPDE}
\dbL^{t,\o, \pi_n} v_n (\pi_n; s, x) := \pa_t v_n(\pi_n; s,x) + G(s, \o\otimes_t \o^{\pi_n}, v_n, \pa_x v_n, \pa^2_{xx} v_n) =0,\q (s, x) \in Q^\e_{t_n},
\eea
where $\o^{\pi_n}$ is the linear interpolation of $(t, {\bf 0}), (t_i, \sum_{j=1}^i x_j)_{1\le i\le n}, (T, \sum_{j=1}^n x_j)$,

(ii) $u^{t, \o, \e} \le u^{t,\o}$ on $\L^t$ and $\lim_{\e\to 0} u^{t, \o, \e}(t, {\bf 0}) = u(t,\o)$.

We define  pseudo Markovian viscosity super-solution similarly , and we call $u$ a pseudo Markovian viscosity solution if it is both a  pseudo Markovian viscosity sub-solution and super-solution.   
\end{defn}

\begin{rem}
\label{rem-viscosity}
{\rm  (i) Definition \ref{defn-viscosity} is completely local. Definition  \ref{defn-viscosity2} is in between local and global. The $u^{t, \o, \e}$ may depend on $(t,\o)$ and we require the convergence of $u^{t,\o,\e}$ only at $(t,\o)$.  In this sense our definition is local. However, the viscosity property of $u^{t,\o,\e}$ and the inequality $u^{t, \o, \e} \le u^{t,\o}$ hold on $\L^t$ and in this sense the definition is global. 

(ii) In Definition \ref{defn-viscosity2} the  $v_n$ is required to satisfy the path frozen PDE \reff{vnPDE}, and thus this definition relies heavily on the path frozen PDE. We will see in Proposition \ref{prop-equivalence} below that, when $G$ is uniformly continuous in $\o$,  one can give an equivalent definition using the original PPDE \reff{PPDE}. 

(iii) Both \cite{ETZ2} and \cite{RTZ3} require the uniform continuity of $G$ in $\o$, which is not required in this paper. We remark that this uniform regularity can be violated even in semilinear case: $G = {1\over 2} \si^2(t,\o) : \g + f(t,\o, y,z)$. 
\qed}
\end{rem}

In the following proposition we state an alternative definition which is equivalent to Definiton \ref{defn-viscosity2} when the generator $G$ is uniformly continuous in $\o$.  The proof is postponed to Appendix.
 \begin{prop}
 \label{prop-equivalence}
Let Assumption \ref{assum-G} hold true, and assume further that $G$ is uniformly continuous in $\o$. Then $u$ is a pseudo Markovian $\cP_L$-viscosity sub-solution if and only if, for any $(t,\o)\in [0,T)\times \O$, there exist $u^{t, \o, \e} \in \cM_\e(\L^t)$, $\e>0$, such that 

(i)  for all $(t',\o')\in [t, T)\times \O^t$, $u^{t, \o, \e}$ is a viscosity sub-solution of PPDE \reff{PPDE} at $(t', \o\otimes_t \o')$;

(ii) $u^{t, \o, \e} \le u^{t,\o}$ on $\L^t$ and $\lim_{\e\to 0} u^{t, \o, \e}(t, {\bf 0}) = u(t,\o)$.
 \end{prop}

 \section{Comparison principle}
\label{sect-comparison}
\setcounter{equation}{0}
The main result of this paper is the following comparison  principle for pseudo Markovian Crandall-Lions viscosity solutions. Since  a $\cP_L$-viscosity semi-solution is always a Crandall-Lions viscosity semi-solution, it also implies  the comparison principle for pseudo Markovian $\cP_L$-viscosity solutions.
 \begin{thm}
\label{thm-comparison}
Let Assumption \ref{assum-G} hold. Assume $u_1$ and $u_2$ are  pseudo Markovian Crandall-Lions viscosity sub-solution and super-solution of PPDE \reff{PPDE}, respectively. If  $u_1(T,\cd)\le u_2(T,\cd)$, then $u_1\le u_2$ on $\L$.
\end{thm}
\proof In this proof, viscosity semi-solutions are always in  Crandall-Lions sense. 
Without loss of generality, we shall only prove $u_1(0,{\bf 0}) \le u_2(0,{\bf 0})$. For $i=1,2$, let $u_i^\e\in \cM_\e(\L)$ be the corresponding approximations with corresponding $v^i_n$.  By Definition \ref{defn-viscosity2} (ii),  it suffices to show that $u^\e_1(0,{\bf 0}) \le u^\e_2(0,{\bf 0})$ for all $\e>0$.  In the rest of this proof we fix $\e>0$ and denote 
$w_n := v^1_n-v^2_n$.

{\it Step 1.} We first show that, for any $n\ge 0$, $\pi_n = (t_i, x_i)_{1\le i\le n} \in \Pi^\e_n$,  it holds 
\bea
\label{wn}
w^+_n(\pi_n; t_n, {\bf 0})\le \overline \cE^L_{t_n} \left[e^{L (\ch-t_n)}w^+_n(\pi_n; \ch, B^{t_n}_\ch)\right]\q\mbox{where}\q \ch := \ch^{t_n,0, \e}_1.
\eea
Without loss of generality it is enough to prove the statement for $n=0$. That is, denoting $w:=w_0$, 
\bea
\label{w}
w^+(0, {\bf 0})\le \overline \cE^L \left[e^{L \ch}w^+( \ch, B_\ch)\right]\q\mbox{where}\q \ch := \ch^{\e}_1.
\eea

For any $\d>0$ small, by Remark \ref{rem-Markov} (ii), $v^1_0, v^2_0$ are uniformly continuous viscosity semi-solutions of  the following PDE: 
\beaa
\pa_t v + G(t, {\bf 0}, v, \pa_x v, \pa^2_{xx} v) =0,\q (t,x)\in D^{\e,\d}_1.
\eeaa
Following the arguments in \cite{ETZ2}, Lemma 6.1, or following an alternative argument in  \cite{BBP}, Lemma 3.7, one can easily prove that
for any $(t,x) \in D^{\e,\d}_1$,
\beaa
w^+(t,x) \le   \overline \cE^L_t \Big[ e^{L (\ch^{t, x,\e-L_1t}_1 - t)}w^+\big(\ch^{t, x,\e-L_1t}_1,B^t_{\ch^{t, x,\e-L_1t}_1}\big)\Big]. 
\eeaa
In particular, this implies
\beaa
w^+(\d,\o_\d)\le   \overline\cE^L_\d \Big[e^{L (\ch^{\d, x,\e-L_1\d}_1 - \d)}w^+\big(\ch^{\d, \o_\d,\e-L_1\d}_1,B^\d_{\ch^{\d, \o_\d,\e-L_1\d}_1}\big)\Big],\q\mbox{if}~ \ch(\o) >\d. 
\eeaa
By the uniform regularity of $w$, and using \eqref{chDPP} we see that
\beaa
\ol\cE^L\Big[w^+(\d,B_\d) \1_{\{\ch >\d\}}\Big] \le \overline \cE^L \Big[e^{L (\ch - \d)}w^+\big(\ch,B_\ch\big)\1_{\{\ch >\d\}}\Big] \le \overline \cE^L \Big[e^{L \ch}w^+\big(\ch,B_\ch\big)\Big] .
\eeaa
Then, noting that $|w|\le C$, 
\beaa
&&w^+(0,0) -  \overline \cE^L \Big[e^{L\ch} w^+\big(\ch,B_\ch\big)\Big] \le w^+(0,0)- \ol\cE^L\Big[w^+(\d,B_\d) \1_{\{\ch >\d\}}\Big]  \\
&\le& \ol\cE^L\Big[|w^+(\d,B_\d)- w^+(0,0)|\Big] + C\overline \cE^L \big[\1_{\{\ch >\d\}}\big] \\
&\le&  \ol\cE^L\Big[|w^+(\d,B_\d)- w^+(0,0)|\1_{\{|B_\d|\le \d^{1\over 3}\}}\Big] + C\ol\cE^L\big[\1_{\{|B_\d|> \d^{1\over 3}\}}\big]+ C\overline \cE^L \big[\1_{\{\ch >\d\}}\big] \to 0,\q\mbox{as}~\d\to 0,
\eeaa
where the first convergence is due to the continuity of $z$ at $(0,0)$, the second one is due to standard estimates, and the third one is due to \reff{chlowerbound}.

{\it Step 2.} We next show that, for any $n\ge 0$, 
\bea
\label{Ewn}
\ol\cE^L\Big[e^{L\ch^\e_n}w^+_n(\pi^\e_n(B); \ch^\e_n, {\bf 0})\Big]\le \overline \cE^L \left[e^{L \ch^\e_{n+1}}w^+_{n+1}(\pi^\e_{n+1}(B); \ch^\e_{n+1}, 0)\right].
\eea
Indeed, for any $\o$, by \reff{wn} and \reff{compatibility} we have:  denoting $\ch^{\o,\e}_{n+1} := (\ch^{\e}_{n+1})^{\ch^\e_n(\o), \o}$,
\beaa
e^{L\ch^\e_n(\o)}w^+_n(\pi^\e_n(\o); \ch^\e_n(\o), {\bf 0})\le \overline \cE^L_{\ch^\e_n(\o)} \left[e^{L \ch^{\o,\e}_{n+1}}w^+_{n+1}\big(\pi^\e_{n}(\o), (\ch^{\o,\e}_{n+1}, B^{\ch^\e_n}_{\ch^{\o,\e}_{n+1}}) ; \ch^{\o,\e}_{n+1}, 0\big)\right], 
\eeaa
Now for any $\d>0$ and $\dbP\in \cP$, by the uniform regularity in Remark \ref{rem-Markov} (ii), it follows from the arguments in  \cite{ETZ2} (5.5) that, 
\beaa
&&e^{L\ch^\e_n}w^+_n(\pi^\e_n; \ch^\e_n, {\bf 0})\1_{\cap_{i=1}^n\{\ch^\e_{i}-\ch_{i-1}^\e\ge \d\}}\\
&\le& \esssup^\dbP_{\dbP'\in\cP_L (\dbP,\ch^\e_n)}\dbE^{\dbP'} \left[e^{L \ch^{\e}_{n+1}}w^+_{n+1}\big(\pi^\e_{n+1}(B); \ch^{\e}_{n+1}, 0\big)\1_{\cap_{i=1}^n\{\ch^\e_{i}-\ch_{i-1}^\e\ge \d\}}\Big|\cF_{\ch^\e_n}\right], \q\dbP\mbox{-a.s.}
\eeaa
where $\cP_L(\dbP,\ch^\e_n):= \{\dbP'\in \cP_L:\dbP=\dbP' \mbox{ on }\cF_{\ch^\e_n}\}$. This implies
\beaa
&&\ol\cE^L\Big[e^{L\ch^\e_n}w^+_n(\pi^\e_n(B); \ch^\e_n, {\bf 0}) \1_{\cap_{i=1}^n\{\ch^\e_{i}-\ch_{i-1}^\e\ge \d\}} \Big]\\
&\le& \overline \cE^L \left[e^{L \ch^\e_{n+1}}w^+_{n+1}(\pi^\e_{n+1}(B); \ch^\e_{n+1}, 0) \1_{\cap_{i=1}^n\{\ch^\e_{i}-\ch_{i-1}^\e\ge \d\}}\right].
\eeaa
Recall \reff{chlowerbound} and send $\d\to 0$, 
we obtain \reff{Ewn} immediately.

{\it Step 3.} Applying Step 2 repeatedly, we have
\bea
\label{west}
w^+(0,0) \le \overline \cE^L \left[e^{L \ch^\e_{n}}w^+_{n}(\pi^\e_{n}(B); \ch^\e_{n}, 0)\right]  \le e^{LT}\overline \cE^L \left[w^+_{n}(\pi^\e_{n}(B); \ch^\e_{n}, 0)\right],\q\forall n\ge 0.
\eea
Note that $w^+_{n}(\pi^\e_{n}(B); \ch^\e_{n}, 0)=0$ on $\{\ch^\e_n = T\}$. Send $n\to\infty$ and apply Lemma \ref{lem-chn}, we have
\beaa
w^+(0,0) \le C e^{LT}\overline \cE^L \left[\1_{\{\ch^\e_n < T\}}\right]\to 0.
\eeaa
This completes the proof.
\qed

\section{Existence}
\label{Existence}
\setcounter{equation}{0}

In this section we provide a generic existence result.  

\subsection{A bounding equation}
\label{sect-bounding}
We first investigate a bounding equation which will be crucial for our existence result:
\bea
\label{dbLbar}
\left.\ba{c}
\dis \ol\dbL \ol v(t,x) := \pa_t \overline v(t,x)+\overline g(\ol v,\pa_x \ol v,\pa_{xx}\ol v)=0,\\
\dis \overline g(y,z,\g):={1\over 2} \sup_{ 0\leq \si \le \sqrt{2L}I_d}[\si^2:\g]+L [|y|+|z|]+C_0,
\ea\right.
\eea
where the control $\si$ takes values in $\dbS^d$.  We start with a local result:

\begin{lem}
\label{lem-barg}
Let $t_0<T$, $\e>0$, and $h\in C^0_b(\pa Q^\e_{t_0})$. Then the PDE
\beaa
\ol\dbL \ol v(t,x) =0, ~(t,x)\in Q^\e_{t_0};\q \ol v(t,x) = h(t,x),~ (t,x)\in \pa Q^\e_{t_0},
\eeaa
has a $\cP_L$-viscosity solution $\ol v \in C^0_b(\wh Q^\e_{t_0})$: 
\bea
\label{barv}
&\dis\ol v(t,x) := \sup_{b\in \dbL^0_L(\dbF^{t})} \ol \cE_t^{L}\Big[e^{\int_t^{\ch^{t,x}}b_rdr} h(\ch^{t,x},x+B^t_{\ch^{t,x}})+C_0\int_t^{\ch^{t,x}}e^{\int_s^{\ch^{t,x}}b_rdr} ds\Big],&\\
&\dis\mbox{where}~ \ch^{t,x} := \ch^{t,x, \e-L_1(t-t_0)},\q  {\dbL_L^0(\dbF^t):=\{b\in \dbL^0(\dbF^{t}): |b|\le L\}}.&\nonumber
\eea
\end{lem}
\proof Assume $t_0=0$ for simplicity. First, it is clear that $|\ol v|\le C$.  For any $0<\d<{\e\wedge T\over L_1}$, denote 
\beaa
\ol Q^{\e,\d}_{0} := \Big\{(t,x) \in \wh Q^\e_{t_0}: t\ge \d\Big\}.
\eeaa
Note that  $\pa Q^\e_{0} \cap \ol Q^{\e,\d}_{0} $ is compact, then $h$ is  uniformly continuous on it. By  Lemma \ref{lem-cHLip} one may easily show that the $\ol v$ defined by \reff{barv} is (uniformly) continuous in $\ol Q^{\e,\d}_{0}$. Provided this  regularity, it follows from standard arguments that $\ol v$ satisfies the dynamical programming principle, which implies  further  that $\ol v$ is a viscosity solution of PDE \reff{dbLbar} in $Q^\e_{0} \cap \ol Q^{\e,\d}_{0}$. By the arbitrariness  of $\d$, we see that $\ol v$ is continuous in $Q^\e_{0} \cap \pa Q^{\e}_{0}$ and is a  viscosity solution of PDE \reff{dbLbar} in $Q^\e_{0}$.

It remains to prove the continuity at $(0, 0)$. Let $(t,x) \in Q^\e_{0}$ and denote $\d := t + |x|$. Let
\beaa
\ch:= \ch^{0, 0,\e},\q \tilde \ch := \inf\big\{s\ge t: |x+ B_{s}-B_t| + L_1(s-t) \ge \e-L_1t\big\}\wedge T.
\eeaa
Then one can easily see that
\beaa
\ol v(0,0) &:=& \sup_{b\in \dbL^0_L(\dbF)} \ol \cE^{L}\Big[e^{\int_0^{\ch}b_rdr} h(\ch,B_{\ch})+C_0\int_0^{\ch}e^{\int_s^{\ch}b_rdr} ds\Big];\\
\ol v(t,x) &:=& \sup_{b\in \dbL^0_L(\dbF)} \ol \cE^{L}\Big[e^{\int_t^{\tilde \ch}b_rdr} h(\tilde \ch,x+B_{t,\tilde \ch})+C_0\int_t^{\tilde \ch}e^{\int_s^{\tilde \ch}b_rdr} ds\Big].
\eeaa
Let $\d_0 >0$ be a constant which will be specified later and assume $\d\le \d_0$. Denote 
\beaa
\ol v'(0,0) &:=& \sup_{b\in \dbL^0_L(\dbF)} \ol \cE^{L}\Big[e^{\int_0^{\ch}b_rdr} h(\ch\vee \d_0,B_{\ch})+C_0\int_0^{\ch}e^{\int_s^{\ch}b_rdr} ds\Big];\\
\ol v'(t,x) &:=& \sup_{b\in \dbL^0_L(\dbF)} \ol \cE^{L}\Big[e^{\int_t^{\tilde \ch}b_rdr} h(\tilde \ch\vee \d_0,x+B_{t,\tilde \ch})+C_0\int_t^{\tilde \ch}e^{\int_s^{\tilde \ch}b_rdr} ds\Big].
\eeaa
Here $h(\d_0, x) := h(\d_0, {\e-L_1\d_0\over |x|}x)$ when $|x|> \e-L_1\d_0$.    Then, for $\d \le {\e\over 2}\wedge \d_0$, by \reff{chlowerbound} we have 
\beaa
|\ol v'(0,0)-\ol v(0,0)| \le C\sup_{\dbP\in \cP^L} \dbP(\ch \le \d_0)\le C_\e\d_0;\q |\ol v'(t, x)-\ol v(t, x)| \le C\sup_{\dbP\in \cP^L} \dbP(\tilde \ch \le \d_0) \le C_\e \d_0.
\eeaa
Notice that $h$ is uniformly continuous on $\ol Q^{\e,\d_0}_{0}$, by Lemma \ref{lem-cHLip} one can easily prove that
\beaa
|\ol v'(0,0)-\ol v'(t,x)| \le C\rho_{\d_0}(\d),
\eeaa
where the modulus of continuity function $\rho_{\d_0}$ may depend on $\d_0$.  Thus, for any $\d_0>0$,
\beaa
|\ol v(0,0)-\ol v(t,x)| \le C\rho_{\d_0}(\d) + C_\e \d_0.
\eeaa
This implies  $\dis\limsup_{\d\to 0} |\ol v(0,0)-\ol v(t,x)| \le C_\e \d_0$. Since $\d_0$ is arbitrary,  we obtain $\dis\lim_{\d\to 0} \ol v(t,x) = \ol v(0, 0)$.
\qed

We next extend the above construction to a global one on $[0, T]$.  Our construction is similar to that in \cite{ETZ1}, Section 7.  Given $(\pi_n, t, x) \in D^\e_{n+1}$, define
\bea
\label{existence-chn0}
\ch^{t,x}_1 := \ch^{t,x, \e-  L_1(t-t_n)},\q \ch^{t,x}_{m+1} := \ch^{\ch_m^{t,x}, 0, \e}(B^t_\cd- B^t_{\ch^{t,x}_m}),~m\ge 1, 
\eea
and let $B^{\e, \pi_n,t,x}$ be the linear interpolation of: denoting $t_0:=0$,
\bea
\label{Be0}
(t_i, \sum_{j=1}^ix_j)_{0\le i\le n}, \q   (\ch^{t,x}_m,  \sum_{i=1}^n x_i+x+ B^t_{\ch^{t,x}_m})_{m\ge 1}.
\eea
We then define
\bea
\label{olthen0}
 \ol\th^\e_n(\pi_n; t,x) :=\sup_{b\in \dbL^0_L(\dbF^t)} \ol \cE^{L}\Big[e^{\int_t^T b_rdr} \xi(B^{\e,\pi_n,t,x}) +C_0\int_t^Te^{\int_s^T b_rdr} ds\Big].
\eea
Our main result of this subsection is:

\begin{lem}
\label{lem-olthen}
For any $\xi\in UC_b(\O)$, $\ol \th^\e_n$ is bounded and continuous in $D^\e_{n+1}$.
\end{lem}

We remark that in general $\ol\th^\e_n$ may not be continuous on the closure of $D^\e_{n+1}$. The proof of Lemma \ref{lem-olthen} relies heavily on the regularity of the hitting times established in Lemma \ref{lem-cHLip}. It is quite lengthy and we postpone it to Appendix.

\subsection{A general existence result}
We shall assume
\begin{assum}
\label{assum-existence}
For any $\e>0$, $\pi_n\in \Pi^\e_n$,  and $h\in C^0_b(\pa Q^\e_{t_n})$, the PDE \reff{vnPDE} with boundary condition $h$ has a $\cP_L$-viscosity solution (resp. Crandall-Lions viscosity solution) $v_n \in C^0_b(\wh Q^\e_{t_n} )$.
\end{assum}

The following lemma is crucial.
\begin{lem}
\label{lem-then}
 Let Assumptions \ref{assum-G}  and \ref{assum-existence} hold.  Then there exists $u^\e \in \cM_\e(\L)$  with corresponding $\th^\e_n := v_n$ such that

(i) $\th^\e_n(\pi_n; \cd)$ is a $\cP_L$-viscosity solution (resp. Crandall-Lions viscosity solution) to PDE \reff{vnPDE};

(ii) $u^\e(T, \o) =  \xi(\o^{(\ch^\e_n, \o_{\ch^\e_n}-\o_{\ch^\e_{n-1}})_{n\ge 1}})$.   
 \end{lem}
\proof We follow the  arguments in \cite{ETZ2}, Lemma 6.3. By Lemmas \ref{lem-barg} and \ref{lem-olthen} we see that $\ol\th^\e_n(\pi_n;\cd)$ is a viscosity solution of PDE \reff{dbLbar} in $Q^\e_{t_n} \cup \pa Q^\e_{t_n}$, for any given $\pi_n\in \Pi^\e_n$. Introduce similarly 
\bea
\label{ulthen}
 \ul\th^\e_n(\pi_n; t,x) :=\inf_{b\in \dbL^0_L(\dbF)} \ul \cE^{L}\Big[e^{\int_t^T b_rdr} \xi(B^{\e,\pi_n,t,x}) -C_0\int_t^Te^{\int_s^T b_rdr} ds\Big],
\eea
which corresponds to the lower bounding equation: 
\bea
\label{dbLul}
\left.\ba{c}
\dis \ul\dbL \ul v(t,x) := \pa_t \overline v(t,x)+ \ul g(\ul v,\pa_x \ul v,\pa_{xx}\ul v)=0,\\
\dis \ul g(y,z,\g):={1\over 2} \inf_{ 0\leq \si \le \sqrt{2L}I_d}[\si^2:\g]-L [|y|+|z|]-C_0,
\ea\right.
\eea

For each $m\ge 1$, define two functions on $D^\e_{m+1}$: 
\beaa
\overline \theta_m^{\e,m}(\pi_m;t,x):=\overline \theta_m^{\e}(\pi_m;t,x),\q  \underline \theta_m^{\e,m}(\pi_m;t,x):=\underline \theta_m^{\e}(\pi_m;t,x).
\eeaa
We next define $\ol \theta_i^{\e,m}(\pi_i;\cdot)\in C^0(\wh Q^\e_{t_i})$ and   $\ul \theta_i^{\e,m}(\pi_i;\cdot)\in C^0(\wh Q^\e_{t_i})$, $i=m-1,\cds, 0$, backwardly as the unique viscosity solutions of the following PDE:
\bea
\label{existence-Lv}
\pa_t v(s,x) + G(s,  \o^{\pi_i}, v, \pa_x v, \pa^2_{xx} v) =0,\q (s, x) \in Q^\e_{t_i},
\eea
with boundary conditions $\overline \theta_{i+1}^{\e,m}(\pi_i, {(t,x)};t,0)$  and $\underline \theta_{i+1}^{\e,m}(\pi_i, {(t,x)};t,0)$, respectively. Here the existence of such viscosity solution is due to Assumption  \ref{assum-existence}, while the comparison principle and the uniqueness of viscosity solutions are implied by Assumption \ref{assum-G}. Clearly $\ul \th^\e_{m+1} \le \ol \th^\e_{m+1}$. Since $\ul g \le G(t,\o,\cd) \le \ol g$, by the comparison principle for the PDE \reff{existence-Lv} we see that
\beaa
\ul \th^{\e,m}_m \le \ul \th^{\e, m+1}_m \le  \ol \th^{\e, m+1}_m \le \ol \th_m^{\e, m}.
\eeaa
Then,  applying repeatedly the comparison principle for the PDE \reff{existence-Lv}  on $Q^\e_{t_i}$, backwardly in $i$, 
\bea
\label{monotonicity}
\ul\th^\e_i \le \ul \th^{\e,m}_i \le \ul \th^{\e,m+1}_i  \le  \ol \th^{\e,m+1}_i   \le \ol \th^{\e,m}_i \le \ol \th^\e_i,\q 1\le i\le m.
\eea

Denote $w^{\e,m}_i:=\overline \theta_i^{\e,m} -\underline \theta_i^{\e,m}$ and recall the notations in  \reff{existence-chn0}. Since both $\overline \theta_i^{\e,m}$ and  $\underline \theta_i^{\e,m}$ are viscosity solutions of PDE \reff{existence-Lv} (with different boundary conditions), it follows from the proof of Theorem \ref{thm-comparison}, in particular \reff{west}, that
\beaa
0 \le w^{\e,m}_i(\pi_i; t,x) \le C \ol\cE^L_t\Big[w^{\e,m}_m\big(\pi_i, (\ch^{t,x}_1, x+ B^t_{\ch_1}), (\ch^{t,x}_i, B^t_{\ch^{t,x}_i} - B^t_{\ch^{t,x}_{i-1}})_{2\le i\le m-i}; \ch^{t,x}_{m-i}, 0\big)\Big].
\eeaa
It is clear that $w^{\e,m}_m(\pi_m; t_m,0)=0$ when $t_m=T$. Then, by Lemma \ref{lem-chn},
 \beaa
0 \le w^{\e,m}_i(\pi_i; t,x) \le C \sup_{\dbP\in \cP^t_L} \dbP\Big(\ch^{t,x}_{m-i}<T\Big) \le \frac{C}{(m-i)^2\e^2}\to 0,\q\mbox{as}~m\to\infty.
 \eeaa
Together with the monotonicity (in $m$) in \reff{monotonicity}, this implies the following limits exist and are equal: 
\beaa
\theta_n^\e :=\lim_{m\to\infty} \overline \theta_n^{\e,m}=\lim_{m\to\infty} \underline \theta_n^{\e,m},
\eeaa
  where the first limit is decreasing and the second one increasing. Since $\ol \th^{\e,m}_n$ and $\ul \th^{\e,m}_n$ are continuous in $\wh D^\e_{n+1}$, then as their monotone limit $\theta_n^\e$ is both upper semicontinuous and lower semicontinuous, and consequently $\theta_n^\e$ is also continuous in $\wh D^\e_{n+1}$.  The viscosity property of $\th^\e_n$ follows from the standard stability result, and all other properties can be verified straightforwardly.  
  \qed

Our main existence result is as follows. 
\begin{thm}
\label{thm-existence}
Let Assumptions \ref{assum-G}  and \ref{assum-existence} hold, and $u^\e$ be as in Lemma \ref{lem-then}. If $u^\e$ converges to $u$ uniformly,  then $u$ is a pseudo Markovian $\cP_L$-viscosity solution (resp. Crandall-Lions viscosity solution) of PPDE  \reff{PPDE} with terminal condition $\xi$. 
\end{thm}
\proof Without loss of generality, we may assume $G$ is non-increasing in $y$.  Denote  
\beaa
\rho(\e) := \|u^\e-u\|_\infty,~ \th^{\e,+}_n := \th^\e_n + \rho(\e),~ u^{\e,+} := u^\e + \rho(\e),~  \th^{\e,-}_n := \th^\e_n - \rho(\e),~ u^{\e,-} := u^\e - \rho(\e).
\eeaa
Then $u^{\e,-} \le u\le u^{\e,+}$, and $\lim_{\e\to 0} u^{\e,+} = u =  \lim_{\e\to 0} u^{\e,-}$. It is straightforward to show that $ \th^{\e,+}_n(\pi_n;\cd)$ and $ \th^{\e,-}_n(\pi_n;\cd)$ are viscosity super solution and viscosity subsolution of PDE \reff{vnPDE}, respectively.  This implies that $u$ is a   pseudo Markovian viscosity solution of PPDE  \reff{PPDE}. Finally, it is obvious that  $u(T,\cd) = \xi$. 
\qed

\section{Stochastic HJB equations} 
\label{sec-HJB}
\setcounter{equation}{0}
The stochastic HJB equation is introduced in \cite{Peng-SHJB}  to characterize the value function for an optimization problem with random coefficients. Let $\dbU$ be an arbitrary measurable set, $\cU^t$ be the set of $\dbF^t$-progressively measurable and $\dbU$-valued processes. Given $(t,\o, x)\in [0, T)\times \O\times \dbR^{d'}$ for some dimension $d'$  and $\a\in \cU^t$, consider the following controlled decoupled FBSDE: 
\bea
\label{HJB-FBSDE}
\left.\ba{lll}
\dis X_s=x+\int_t^s b^{t,\o}(r,B^t, X_r,\a_r)dr +\int_t^s \si^{t,\o}(r,B^t,X_r,\a_r)dB^t_r;\\
\dis Y_s=g^{t,\o}(B^t,X_T)+\int_s^T f^{t,\o}(r,B^t,X_r,Y_r,Z_r,\a_r)dr -\int_s^T Z_rdB_r;
\ea\right. t\le s\le T, ~\dbP_0^t\mbox{-a.s.}
\eea
Here $Y$ is a scalar process, $b, \si, f, g$ have appropriate dimensions, $b, \si, f$ are $\dbF$-progressively measurable in all variables, and $g$ is $\cF_T\times \cB(\dbR^{d'})$-measurable.   We shall assume

\begin{assum}
\label{assum-HJB}
(i) $b(t, \o, x, \a)$, $\si(t, \o, x, \a)$, $f(t, \o, x,0, 0,\a)$, and $g(\o,x)$ are bounded;

(ii) $b$ and $\si$ are uniformly Lipschitz continuous in $x$, and $f$ is uniformly Lipschitz continuous in $(y,z)$;

(iii) $b$ and $\si$ are uniformly continuous in $\o$;  and $f$ and $g$ are uniformly continuous in $(\o, x)$; 

(iv) $b, \si$ and $f$ are continuous in $t$.
\end{assum}

Under the above conditions, it is clear that the decoupled FBSDE \reff{HJB-FBSDE} has a unique solution, denoted as $(X^{t,\o,x,\a}, Y^{t,\o,x,\a}, Z^{t,\o,x,\a})$.  We now introduce the optimization problem:
\bea
\label{HJB-u}
u^0(t,\o,x):=\sup_{\a\in \cU^t} Y_t^{t,\o,x,\a}.
\eea
To characterize the above random field $u^0$, \cite{Peng-SHJB} introduced the following stochastic HJB equation (in some simpler case) with $\dbF$-progressively measurable solution pair $(u, q)$:
\bea
\label{HJB-BSPDE}
d u(t,\o,x) &=& -\sup_{\a\in \dbU}\Big[ \frac{1}{2} \pa_{xx}u :\si\si^\top (t,\o,x,\a)+\pa_{x}q : \si(t,\o,x,\a)\nonumber\\
&&+ \pa_x u b(t,\o,x,\a)+f(t,\o,x,u,q+\pa_x u \si(t,\o,x,\a), \a)\Big] + q(t,\o,x) dB_t;\\
u(T,\o,x) &=& g(\o,x).\nonumber
\eea
This is a special type of backward SPDE. When $\si$ does not depend on $\a$, namely there is no diffusion control,  \cite{Peng-SHJB} established its wellposedness in Sobolev sense. The recent work \cite{Qiu} extended the result to the general case with diffusion control, also in terms of Sobolev solutions. 

We shall view the equation \reff{HJB-BSPDE} as a PPDE, as  shown in \cite{ETZ1}, Section 6:
\bea
\label{HJB-PPDE}
&\cL u(t,\o,x) =0,\q u(T,\o,x) = g(\o, x),\q\mbox{where}&\\
&\cL u(t,\o,x) := \pa_t u + \frac{1}{2}\tr(\pa_{\o\o}u)+\sup_{\a\in \dbU}\Big[ \frac{1}{2} \pa_{xx}u :\si^2 (t,\o,x,\a)+\pa_{x\o}u : \si(t,\o,x,\a)&\nonumber\\
& + \pa_x u b(t,\o,x,\a)+f(t,\o,x,u,\pa_\o u+\pa_x u \si(t,\o,x,\a), \a)\Big].&\nonumber
\eea
Indeed, if $u$ is smooth, by comparing \reff{HJB-BSPDE} and the functional It\^{o} formula \reff{Ito} and noting that $d\la B\ra_t = I_d dt$, $\dbP_0$-a.s., one may obtain $q = \pa_\o u$ and then \reff{HJB-PPDE} follows immediately.  In general, of course, $u^0$ is not smooth. Our goal is to characterize $u^0$ as the unique pseudo Markovian viscosity solution of \reff{HJB-PPDE}.  We remark that, the Sobolev theory in \cite{Peng-SHJB, Qiu} requires the special structure of HJB equation, and may not be easy to extend to more general cases like Isaacs equations induced from games. Our viscosity solution characterization, however, can be easily extended as we will see in next section.   

The PPDE  \reff{HJB-BSPDE} is slightly different from \reff{PPDE} due to the involvement of the additional variable $x$. In \cite{ETZ1} we  view $x$ as the current value of another path $\tilde \o$, namely we increase the dimension of the canonical space with canonical paths $(\o, \tilde \o)$, and consider only viscosity solutions in the form $u(t, \o_\cd, \tilde \o_t)$. We emphasize that this PPDE is always degenerate, and thus is not covered by the comparison result in \cite{ETZ2}. The results in this paper apply to this case, and $u^0$ is indeed the unique pseudo Markovian viscosity solution of PPDE \reff{HJB-PPDE}.

However, in this subsection we shall treat \reff{HJB-PPDE} in an alternative way. Note that the main feature of our new notion is the Markovian structure. Since the $x$ part is already Markovian, there is no need to introduce an additional path $\tilde \o$ and then discretize  it. So we shall discretize $\o$ only and deal with $x$ directly. For notational clarity, we will use $\bar x$ to denote the variable induced from the discretization of  $\o$.  Recall the notations in \reff{che} and \reff{Oet} corresponding to the discretization of $\o$. 
Analogous to Definitions \ref{defn-Markov} and \ref{defn-viscosity2}, we define:
\begin{defn}
\label{defn-HJB-Markovian}
Let $\e>0$. We say  $u\in \dbL^0(\L\times \dbR^{d'})$ is $\e$-Markovian, denoted as $u\in \cM_\e(\L\times \dbR^{d'})$, if there exist  deterministic functions $v_n: D^\e_{n+1}\times \dbR^{d'} \to \dbR$, $n\ge 0$, satisfying: 

(i) \reff{Cbar12} holds, namely
\bea
\label{HJB-Markov}
u (t,\o,x)=  \sum_{n=0}^\infty v_n\left(\pi^\e_n(\o); t, \o_t-\o_{\ch^\e_n (\omega)}, x \right) \1_{\{\ch^\e_n(\o)\le t< \ch^\e_{n+1}(\o) ~\mbox{or}~\ch^\e_n(\o) < \ch^\e_{n+1}(\o)=T=t\}}.
\eea

(ii) For all $\pi_n = (t_i, x_i)_{1\le i\le n}\in \Pi^\e_n$ and $(t,\bar x,x)\in \pa Q^\e_{t_n}\times \dbR^{d'}$,  the following compatibility condition holds 
\bea
\label{HJB-compatibility}
v_n(\pi_n; t, \bar x,x)= v_{n+1}(\pi_n, (t,\bar x); t,0,x).
\eea

(iii) Each $v_n$, $n\ge 0$, is continuous in $D^\e_{n+1}\times \dbR^{d'}$.
\end{defn}

\begin{defn}
\label{defn-HJB-viscosity}
We say $u$ is a pseudo Markovian $\cP_L$-viscosity sub-solution (resp. Crandall-Lions viscosity sub-solution)  of PPDE \reff{HJB-PPDE} at $(t,\o,x) \in [0, T)\times \O\times \dbR^{d'}$ if there exist $u^{t, \o, \e} \in \cM_\e(\L^t\times \dbR^{d'})$ with corresponding $\{v_n, n\ge 1\}$,  $\e>0$, such that

(i)  for each $\pi_n = (t_i, x_i)_{1\le i\le n} \in \Pi^{t,\e}_n$, $v_n(\pi_n; \cd)$ is a $\cP_L$-viscosity sub-solution (resp. Crandall-Lions viscosity sub-solution)  to the following PDE: at $(\pi_n, s, \bar x, x) \in   \Pi^{t,\e}_n \times Q^\e_{t_n}\times \dbR^{d'}$,
\bea
\label{HJB-vnPDE}
\dbL^{t,\o, \pi_n} v_n&:=& \pa_t v_n + {1\over 2} \pa^2_{\bar x\bar x} v_n:I_d +\sup_{\a\in \dbU}\Big[ \frac{1}{2} \pa_{xx}u :\si^2 (t,\o\otimes_t \o^{\pi_n},x,\a)+\pa_{x \bar x}u : \si(t,\o\otimes_t \o^{\pi_n},x,\a)\nonumber\\
&&+ \pa_x u b(t,\o\otimes_t \o^{\pi_n},x,\a)+f(t,\o\otimes_t \o^{\pi_n},x,u,\pa_{\bar x} u,\pa_x u, \a)\Big] =0.
\eea
where $\o^{\pi_n}$ is the linear interpolation of $(t, {\bf 0}), (t_i, \sum_{j=1}^i x_j)_{1\le i\le n}, (T, \sum_{j=1}^n x_j)$,

(ii) $u^{t, \o, \e} \le u^{t,\o}$ on $\L^t\times \dbR^{d'}$ and $\lim_{\e\to 0} u^{t, \o, \e}(t, {\bf 0},x) = u(t,\o,x)$.

We define  pseudo Markovian viscosity super-solution similarly , and we call $u$ a pseudo Markovian viscosity solution if it is both a  pseudo Markovian viscosity sub-solution and super-solution.   
\end{defn}

Our main result of this section is:
\begin{thm}
\label{thm-HJB}
Let Assumption \ref{assum-HJB} hold and $L>0$ be large enough. Then the $u^0$ defined by \reff{HJB-u} is the unique pseudo Markovian viscosity solution (both in $\cP_L$-sense and in Crandall-Lions sense)   of PPDE \reff{HJB-PPDE}. 
\end{thm}
\proof  Clearly Assumption \ref{assum-HJB} implies Assumption \ref{assum-G}, with the state space raised to $(\o, x)$ due to the involvement of $x$. Then the comparison principle follows the same arguments as in Theorem \reff{thm-comparison}, which implies the uniqueness immediately.  So it suffices to verify that $u^0$ is indeed a pseudo Markovian viscosity solution. Without loss of generality, we shall only verify the viscosity property at $(0,0)$. We note that, due to the representation \reff{HJB-u}, we shall construct the functions $\th^\e_n$ in Lemma  \ref{lem-then} directly, without referring to Assumption \ref{assum-existence}. We next verify the uniform convergence of the corresponding $u^\e$, and then the existence follows from Theorem \ref{thm-existence}.

Our construction of the functions $\th^\e_n$ is similar to that in Subsection \ref{sect-bounding}. Fix $\e>0$, and let $(\pi_n, t, \bar x, x) \in D^\e_{n+1} \times \dbR^{d'}$. Recall the notations in \reff{existence-chn0} and \reff{Be0}, and denote
\bea
\label{HJB-che}
\ch_0:= t_n,\q \ch_m := \ch^{t, \bar x}_m, ~m\ge 1,\q \wh B := B^{\e, \pi_n,t,\bar x}.
\eea
 We then define
\bea
\label{HJB-then}
\th^\e_n(\pi_n; t, \bar x, x) := \sup_{\a \in \cU^t} Y^{\pi_n, t, \bar x, x, \a}_t,
\eea
where $(X^{\pi_n, t, \bar x, x, \a}, Y^{\pi_n, t, \bar x, x, \a}, Z^{\pi_n, t, \bar x, x, \a})$ is the solution to the following decoupled FBSDE on $[t, T]$:
\bea
\label{HJB-FBSDE2}
\left.\ba{lll}
\dis X_s=x+\int_t^s \sum_{m=0}^\infty \Big[b(r,\wh B_{\cd\wedge \ch_m}, X_r,\a_r)\1_{[\ch_m, \ch_{m+1})}dr + \si(r, \wh B_{\cd\wedge \ch_m},X_r,\a_r)\1_{[\ch_m, \ch_{m+1})}dB^t_r\Big];\\
\dis Y_s=g(\wh B, X_T)+\int_s^T \sum_{m=0}^\infty f(r,\wh B_{\cd\wedge \ch_m},X_r,Y_r,Z_r,\a_r)\1_{[\ch_m, \ch_{m+1})} dr -\int_s^T Z_rdB_r,\q \dbP_0^t\mbox{-a.s.}
\ea\right. 
\eea
While it is not completely trivial, it follows from similar arguments in Lemmas \ref{lem-barg} and  \ref{lem-olthen} that  $\th^\e_n$ satisfies all the requirements in Lemma  \ref{lem-then}.  We leave the details to interested readers.

Moreover,   as in \reff{HJB-Markov} we denote
\bea
\label{HB-ue}
u^\e(t,\o,x) := \sum_{n=0}^\infty \th^\e_n (\pi_n^\e(\o);  t, \o_t - \o_{\ch^\e_n(\o)}, x)  \1_{\{\ch^\e_n(\o)\le t< \ch^\e_{n+1}(\o) ~\mbox{or}~\ch^\e_n(\o) < \ch^\e_{n+1}(\o)=T=t\}}.
\eea
It remains to verify that $u^\e$ converges to $u^0$ uniformly.  Indeed, for any $(t,\o,x)$ and $\e>0$, fix the $n$ such that $\ch^\e_n(\o)\le t < \ch^\e_{n+1}(\o)$. For $\pi_n := \pi_n^\e(\o)$, we have
\beaa
\sup_{t\le r\le T} \|\o \otimes_t B^t - \sum_{m=0}^\infty \wh B_{\cd\wedge \ch_m} \1_{[\ch_m, \ch_{m+1})}(r)\|_r  = \sup_{m\ge 0} \sup_{\ch_m \vee t \le r\le \ch_{m+1}} \|\o \otimes_t B^t - \wh B_{\cd\wedge \ch_m} \|_r \le \e.
\eeaa
Since $b$ and $\si$ are uniformly continuous in $\o$, by standard SDE arguments we have
\beaa
\dbE^{\dbP_0^t} \Big[\|X^{\pi_n, t, \bar x, x, \a} - X^{t,\o,x,\a}\|_T^2 \Big] \le C \rho_1(\e),
\eeaa
for some modulus of continuity function $\rho_1$. Moreover, since $f$ and $g$ are uniformly continuous in $(\o, x)$, by standard BSDE arguments we obtain
\beaa
  \dbE^{\dbP_0^t} \Big[\|Y^{\pi_n, t, \bar x, x, \a} - Y^{t,\o,x,\a}\|_T^2 + \int_t^T |Z^{\pi_n, t, \bar x, x, \a}_r - Z^{t,\o,x,\a}_r|^2dr \Big] \le C \rho_2(\e),
\eeaa
for some modulus of continuity function $\rho_2$. By the arbitrariness of $\a$, this implies that 
\beaa
|u^\e(t,\o,x) - u^0(t,\o,x)|\le C\rho_2(\e).
\eeaa
Now by Theorem \ref{thm-existence} we see that $u^0$ is a pseudo Markovian viscosity solution of PPDE \reff{HJB-PPDE}.
\qed

\section{Path dependent Isaacs equation}
\label{sect-game}
\setcounter{equation}{0}
In this section we study path dependent Isaacs equation, which is the PPDE \reff{PPDE} with  generator:
\bea
\label{gameG}
G(t,\o,y,z,\g) := \inf_{\b\in\dbV}\sup_{\a\in\dbU}\Big[\frac{1}{2}\si\si^\top(t,\o,\a,\b):\g+f(t,\o,y,z\si(t,\o,\a,\b),\a,\b)\Big].
\eea
where $\dbU$ and $\dbV$ are two measurable sets, and $\si$, $f$ are $\dbF$-progressively measurable. We shall assume
 \begin{assum}
\label{assum-game} (i) $\si(t, \o, \a,\b)$, $f(t, \o, 0,0, \a, \b)$, and $\xi$ are bounded;

(ii) $\si$ is unifromly Lipschitz continuous in $\o$, and $f$ is uniformly Lipschitz continuous in $(y,z)$;

(iii) $f$ and $\xi$ are uniformly continuous in $\o$;

(iv) $\si$ and $f$ are continuous in $t$.
\end{assum}

Under  Assumption \ref{assum-game}, clearly $G$ satisfies Assumption \ref{assum-G}. Then it follows from Theorem \ref{thm-comparison} that the path dependent Isaacs equation \reff{PPDE}-\reff{gameG} has at most one pseudo Markovian viscosity solution.  We remark that \cite{PZ2} established the comparison principle for viscosity solutions of this PPDE in the sense of Definition \ref{defn-viscosity}.  However, it followed the approach in \cite{ETZ2} and requires that: (i) $\si$ is uniformly non-degenerate; (ii) $\si$ does not depend on $\o$, and (iii) the dimension $d\le 2$.  None of them is needed  in this paper.

The goal of this section is to construct a pseudo Markovian viscosity solution.
It is well known that Isaacs equation is induced from zero sum stochastic differential games. There are three possible formulations for the game problem: (i) control versus control in strong formulation; (ii) strategy versus control in strong formulation; and (iii) control versus control in weak formulation. We refer to  \cite{PZ2} for detailed discussions on the three formulations. In particular, under the first approach the value function does not satisfy the dynamic programming principle and thus does not provide a representation for the PPDE.  As discussed in \cite{PZ2}, the weak formulation in (iii) has some advantages from practical point of view for games. However, following this approach it is more tricky to obtain the desired regularity of the value function. Since our focus here is not the game problem, but to provide a representation for the solution to PPDE,  we shall use the strong formulation (ii) which is easier for regularity. This approach was initiated by  \cite{FS-game} for PDEs. In a recent work  \cite{Zhang} applied this approach to study the game problem in path dependent case.  Our construction of the pseudo Markovian viscosity solution is in spirit similar to that in \cite{Zhang}. However, \cite{Zhang} uses deterministic time discretization, which is easier for regularities and is sufficient for study of Sobolev type of solutions, but it does not provide uniform pointwise approximations and is not convenient   for our study of viscosity solutions. Indeed, no connection with PPDE is discussed in \cite{Zhang}. 

To this end,  let $\cU^t$, $\cV^t$ denote the sets of $\dbF^t$-measurable $\dbU$-valued, $\dbV$-valued processes, respectively, and $\cB^t$ the set of adapted strategies $\l: \cV^t\to \cU^t$, here adaptedness means:  for any $\t\in \cT^t$,
\bea
\label{strategy}
\mbox{if}~v^1, v^2\in \cV^t~ \mbox{such that}~ v^1_s = v^2_s, ~t\le s\le \t,\q \mbox{then}~\l(v^1)_s = \l(v^2)_s, t\le s\le \t.
\eea
For any $(t,\o)\in \O$ and $(\a,\b)\in \cU^t \times \cV^t$, consider the following path dependent decoupled FBSDEs: 
\bea
\label{game-SDE}
\left.\ba{lll}
X_s &=& \dis \o_t + \int_t^s \si^{t,\o}(r, X_\cd, \a_r, \b_r) dB^t_r,\\
Y_s &=& \dis \xi^{t,\o}(X_\cd) + \int_s^T f^{t,\o}(r, X_\cd, Y_r, Z_r, \a_r, \b_r) dr - \int_s^T Z_r dB^t_r,
\ea\right. t\le s\le T,~ ~\dbP_0^t\mbox{-a.s.}
\eea
Under Assumption \ref{assum-game}, clearly the above FBSDE is wellposed, and we shall denote its unique solution as $(X^{t,\o,\a,\b}, Y^{t,\o,\a,\b}, Z^{t,\o,\a,\b})$. 
We then define
\bea
\label{game-u0}
u^0(t,\o) := \sup_{\l\in \cB^t} \inf_{\b\in \cV^t}  Y^{t,\o,\l(\b),\b}_t.
\eea

\begin{rem}
\label{rem-gamereg}
{\rm In \reff{game-u0}, the controls $(\l(\b), \b)$ do not depend on the variable $\o$, due to its strong formulation. Then, given $(t, \o^i)$, $i=1,2$, we have
\bea
\label{game-ureg}
|u^0(t,\o^1) - u^0(t,\o^2)| \le   \sup_{\l\in \cB^t} \sup_{\b\in \cV^t}  |Y^{t,\o^1,\l(\b),\b}_t- Y^{t,\o^2, \l(\b), \b}_t|,
\eea
and thus the regularity of $u^0$ (in $\o$) follows from standard SDE/BSDE estimates. Under the weak formulation in \cite{PZ2}, the controls $\a, \b$ are feedback type and thus depend on $\o^1, \o^2$. Then we don't have a simple estimate like \reff{game-ureg}, and  the regularity of $u^0$ is indeed  more difficult to establish.
\qed}
\end{rem}

Our main result of this section is:
\begin{thm}
\label{thm-game}
 Under  Assumption \ref{assum-game}, the $u^0$ defined by \reff{game-u0} is the unique pseudo Markovian viscosity solution of PPDE \reff{PPDE}-\reff{gameG} with terminal condition $\xi$.
 \end{thm}
\proof Similar to Theorem \ref{thm-HJB}, it suffices to construct the desired functions $\th^\e_n$ and show that the corresponding process $u^\e$ converges to $u^0$ uniformly.  Again we will only verify the viscosity property at $(0,0)$.

Fix $\e>0$, and let $(\pi_n, t, x) \in D^\e_{n+1}$. Recall \reff{vnPDE} that $\o^{\pi_n}$ denote the linear interpolation of $(0,0)$, $(t_i, \sum_{j=1}^i x_j)_{1\le i\le n}$, and $(T, \sum_{j=1}^n x_j)$. For any $(\a, \b)$, to adapt to the strong formulation, we define $X := X^{\pi_n, t,x,\a,\b}$ and $\ch_m :=  \ch^{\pi_n, t,x,\a,\b}_m$ recursively as follows. First, 
\bea
\label{game-X1}
X^1_s &:=& \sum_{i=1}^n x_i + x + \int_t^s \si(r, \o^{\pi_n}, \a_r, \b_r) dB^t_r, ~t\le s\le T,~\dbP_0^t\mbox{-a.s.}\nonumber\\
\ch_1 &:=& \inf\big\{s\ge t: |X^1_s - \sum_{i=1}^n x_i| + L_1(s-t) \ge \e - L_1(t-t_n)\big\}\wedge T;\\
X_s &:=& X^1_s, ~ t\le s \le \ch_1.\nonumber
\eea
Next, for $m\geq 1$,  
\bea
\label{game-Xm}
\wh X^m &:=& \mbox{linear interpolation of}~ (0, 0), ~(t_i, \sum_{j=1}^i x_j)_{1\le i\le n},~ (\ch_i, X_{\ch_i})_{1\le i\le m}, (T, X_{\ch_m});\nonumber\\
X^{m+1}_s &:=& X_{\ch_m} + \int_{\ch_m}^s \si(r, \wh X^m, \a_r, \b_r) dB^t_r, ~\ch_m\le s\le T,~\dbP_0^t\mbox{-a.s.}\nonumber\\
\ch_{m+1} &:=& \inf\big\{s\ge \ch_m: |X^{m+1}_s - X_{\ch_m}| + L_1(s-\ch_m) \ge \e \big\}\wedge T;\\
X_s &:=& X^{m+1}_s, ~ \ch_m\le s \le \ch_{m+1}.\nonumber
\eea
Denote $\ch_0:= t_n$, $\wh X := \lim_{m\to\infty} X^m$, and let $(Y^{\pi_n,t,x,\a,\b}, Z^{\pi_n,t,x,\a,\b})$ be the unique solution of the following BSDE:
\bea
\label{game-Y}
Y_s = \xi(\wh X) + \int_s^T \sum_{m=0}^\infty f(r,\wh X_{\cd\wedge \ch_m},Y_r,Z_r,\a_r, \b_r)\1_{[\ch_m, \ch_{m+1})} dr -\int_s^T Z_rdB^t_r,\q \dbP_0^t\mbox{-a.s.}
\eea
We then define
\bea
\label{game-then}
\th^\e_n(\pi_n; t,x) := \sup_{\l\in \cB^t} \inf_{\b\in \cV^t} Y^{\pi_n,t,x,\l(\b),\b}_t.
\eea
In the spirit of \reff{game-ureg}, combined with standard SDE/BSDE estimates, one may follow the arguments in Lemma \ref{lem-olthen} to show that $\th^\e_n \in C^0_b(D^\e_{n+1})$. Moreover, provided the above regularity and by standard arguments, see e.g. \cite{FS-game} or \cite{Zhang}, one can prove the dynamic programming principle for $\th^\e_n(\pi_n;\cd)$, which leads to the desired viscosity property immediately. We again leave the details to interested readers.

Finally we prove the convergence of $u^\e_0 = \th^\e_0(0,0)$ with uniform rate. That is, in \reff{game-X1} and \reff{game-Xm} we shall set $n=0$ and $(t,x) = (0,0)$.  Then we have
\bea
\label{game-X}
X_s = \int_0^s   \sum_{m=0}^\infty \si(r,\wh X_{\cd\wedge \ch_m}, \a_r, \b_r)\1_{[\ch_m, \ch_{m+1})}(r) dB_r, ~0\le s\le T,~\dbP_0\mbox{-a.s.}
\eea
and, by the construction of the hitting times $\ch_m$, 
\bea
\label{game-whX}
\|X- \wh X\|_T \le \e.
\eea
Compare \reff{game-X} with the SDE of $X^{0,0,\a,\b}$ in \reff{game-SDE}, it follows from standard SDE arguments that
\bea
\label{game-Xest}
\dbE^{\dbP_0} \Big[\|X- X^{0,0,\a,\b}\|_T^2\Big] \le C \dbE^{\dbP_0}\Big[\|X- \wh X\|_T^2\Big]  \le  C\e^2.
\eea
Moreover, let $\rho$ denote the modulus of continuity function of $f$ and $\xi$ in terms of $\o$. Then, compare \reff{game-Y} with the BSDE of $Y^{0,0,\a,\b}$ in \reff{game-SDE}, it follows from standard BSDE arguments that
\bea
\label{game-Yest}
\dbE^{\dbP_0} \Big[\|Y- Y^{0,0,\a,\b}\|_T^2\Big] &\le& C \dbE^{\dbP_0}\Big[\rho(\|\wh X- X^{0,0,\a,\b}\|_T)^2\Big] \\
&\le& C \dbE^{\dbP_0}\Big[\rho\big(\|\wh X- X\|_T+ \|X- X^{0,0,\a,\b}\|_T\big)^2\Big]\le C\rho'(\e)^2,\nonumber
\eea
for a possibly different modulus of continuity function $\rho'$ which does not depend on the controls $(\a,\b)$. 
This implies that, for the $Y$ corresponding to $(\a, \b) = (\l(\b),\b)$,
\bea
\label{game-conv}
|u^\e_0-u^0_0| \le \sup_{\l\in \cB^0} \sup_{\b\in \cV^0} |Y_0- Y^{0,0,\l(\b),\b}_0| \le C\rho'(\e),
\eea 
which provides the desired convergence of $u^\e_0$ and thus completes the proof.
\qed

We conclude the section with an application on the zero sum game, which is the main result of \cite{Zhang}. Denote by $\cA^0$ the set of adapted strategies $\l: \cU^0\to \cV^0$.
\begin{cor}\label{cor-value-game}
Let Assumption \ref{assum-game} hold and assume further the following Isaacs condition:
\bea
\label{Isaacs}
G(t,\o,y,z,\g) = \tilde G(t,\o,y,z,\g) := \sup_{\a\in\dbU}\inf_{\b\in\dbV}\Big[\frac{1}{2}\si^2(t,\o,\a,\b):\g+f(t,\o,y,z\si(t,\o,\a,\b),\a,\b)\Big].
\eea 
Then the value of the stochastic differential game exists, namely 
\bea
\label{game-value}
u^0_0 = \tilde u^0_0 := \inf_{\l\in \cA} \sup_{\a\in \cU}  Y^{0,0,\a,\l(\a)}_0 .
\eea
\end{cor}
\proof Define  $\tilde u^0(t,\o)$ in the same spirit as $\tilde u^0_0$. Follow the same arguments as in Theorem \ref{thm-game} we see that $\tilde u^0$ is the unique pseudo Markovian viscosity solution of the PPDE \reff{PPDE} with generator $\tilde G$ and terminal condition $\xi$. Since $G = \tilde G$, by the uniqueness of the pseudo Markovian viscosity solution we obtain $u^0 = \tilde u^0$.
\qed

\section{Appendix}
\label{Appendix}
\setcounter{equation}{0}

\subsection{Proof of Proposition \ref{prop-equivalence}}
We shall only prove the equivalence at $(0, {\bf 0})$. The proof for general $(t,\o)$ follows the same argument.   Let $\rho$ denote the modulus of continuity function of $G$ in terms of $\o$. 
Moreover, by the  change variable formula in \cite{ETZ1}, Proposition 3.14, we may assume without loss of generality that
 \bea
 \label{Monotone}
 \mbox{$G$ is nonincreasing in $y$.}
 \eea

 We first assume $u$ satisfies  the sub-solution property  stated at proposition \ref{prop-equivalence} at $(0,{\bf 0})$.  Let $u^\e \in \cM_\e(\L)$ be the approximation given by this property, with corresponding $\{v_n:n\geq 0\}$.  Denote
 \bea
 \label{tildeue}
 \tilde u^\e(t,\o) := u^\e(t,\o) - \rho(\e)[T- t],\q \tilde v_n(\pi_n; t, x) :=  v_n(\pi_n; t, x) - \rho(\e) [T-t].
 \eea
 One may check straightforwardly that
 \beaa
 \mbox{$\tilde u^\e \in M_\e(\L)$ with corresponding $\tilde v_n$, $\tilde u^\e \le u$, and $\lim_{\e\to 0} u^\e(0, {\bf 0}) = u(0, {\bf 0})$.}
 \eeaa
 Then it remains to prove that $\tilde v_n(\pi_n; \cd)$ is a $\cP_L$-viscosity sub-solution of the PDE \reff{vnPDE}, which in this case becomes:
 \bea
\label{vnPDE0}
\dbL^{\pi_n} \tilde v_n (\pi_n; t, x) := \pa_t \tilde v_n(\pi_n; t,x) + G(s, \o^{\pi_n}, \tilde v_n, \pa_x \tilde v_n,  \tilde \pa^2_{xx} v_n) =0,\q (t, x) \in Q^\e_{t_n}.
\eea

To see this, we fix $(t,x) \in Q^\e_{t_n}$.  For any $\tilde \f \in \ul \cA^L \tilde u^\e(t, \o^{\pi_n, (t,x)})$, namely 
\beaa
\tilde \f\in C^{1,2}(\L^t),\q\tilde \f(t,{\bf 0}) - \tilde v_n(\pi_n; t,x) = 0 = \inf_{\t\in\cT^t} \ul \cE^L_t\Big[\tilde \f_{\t\wedge \ch} -   \tilde v_n(\pi_n; \t\wedge \ch, B_{ \t\wedge \ch})\Big],
\eeaa
where we assume without loss of generality that 
\beaa
\ch \le \ch_\e := \inf\{s> t: |x+ B^t_s| + L_1(s-t)\ge \e\}\wedge T.
\eeaa 
Now denote $\f(s, \o) := \tilde \f(s, \o) + \rho(\e) [T-s]$.  Then obviously $\f \in \ul \cA^L u^\e(t, \o^{\pi_n, (t,x)})$, and thus
\beaa
\pa_t \f(t, 0) + G^{t, \o^{\pi_n, (t,x)}}(t, {\bf 0}, \f, \pa_\o\f, \pa^2_{\o\o} \f)\ge  0.
\eeaa
Note that 
\beaa
\pa_t \f = \pa_t \tilde \f - \rho(\e),\q \pa_\o \f = \pa_\o \tilde \f,\q \pa^2_{\o\o} \f = \pa^2_{\o\o} \tilde \f,\q \mbox{and}\q \f(t,{\bf 0}) = \tilde \f(t,{\bf 0}) + \rho(\e)T.
\eeaa
Then, by \reff{Monotone} and evaluating $\f$, $\tilde\f$ and their derivatives at $(t, {\bf 0})$, we obtain
\beaa
 &&\pa_t \tilde \f(t,{\bf 0}) + G(t, \o^{\pi_n}, \tilde \f, \pa_\o \tilde \f,  \tilde \pa^2_{\o\o} \f)\\
 &=& \pa_t  \f(t,{\bf 0})  + \rho(\e)+ G(t, \o^{\pi_n}, \f - \rho(\e)T, \pa_\o \f,   \pa^2_{\o\o} \f)\\
 &\ge& -  G(t, \o^{\pi_n, (t,x)}, \f, \pa_\o\f, \pa^2_{\o\o} \f)+\rho(\e) + G(t, \o^{\pi_n}, \f, \pa_\o \f,   \pa^2_{\o\o} \f) \\
 &\ge&  \rho(\e) - \rho\Big(\| \o^{\pi_n} -  \o^{\pi_n, (t,x)}\|_t\Big) \ge 0.
  \eeaa
 This implies that $\tilde v_n(\pi_n; \cd)$ is an $\cP_L$-viscosity sub-solution of the PDE \reff{vnPDE0}, and thus $u$ is a pseudo Markovian $\cP_L$-viscosity sub-solution at $(0, {\bf 0})$.
 
 Similarly, if $u$ is a  pseudo Markovian $\cP_L$-viscosity sub-solution at $(0, {\bf 0})$ with approximation $u^\e \in \cM_\e(\L)$  in Definition \ref{defn-viscosity2}, one can show that  $u$ satisfies the property at proposition \ref{prop-equivalence} at $(0, {\bf 0})$ with approximation $\tilde u^\e(t,\o) := u^\e(t,\o) - \rho(\e)[T- t]$.
 \qed

\subsection{Proof of Lemma \ref{lem-olthen}}
Our proof here relies heavily on the regularity results in Lemma \ref{lem-cHLip}.  Notice that  the time regularity \reff{chLipt} requires a shift of canonical process.  To facilitate our proof, we extend the canonical space to $\ol \O := \{\o\in C([0, 2T], \dbR^d): \o_0=0\}$, and still denote $B$, $\dbF$, $\cP_L$ etc. in obvious sense.  Given $(\pi_n, t, x) \in D^\e_{n+1}$, define
\bea
\label{existence-chn}
\ch_1 &:=&  \inf\big\{s\ge t:  |x+B_{s-t}| + L_1(s-t) \ge \e- L_1(t-t_n)\big\} \wedge T,\nonumber\\
\ch_{m+1} &:=& \inf\big\{ s \ge \ch_m: |B_{\ch_m+(T-\ch_1), s+(T-\ch_1)}| + L_1 (s-\ch_m) \ge \e\big\}\wedge T,\q~m\ge 1;\\
N &:=&  \inf\{m\ge 1: \ch_m = T\}.\nonumber
\eea
and let $B^\e_\cd(\pi_n,t,x)$ be the path which is the linear interpolation of
\bea
\label{Be}
&(t_i, X_i)_{0\le i\le n}, \q   (\ch_m,  X_{n+m})_{m\ge 1},\q\mbox{where }\ t_0:=0;&\\
&\q   X_i := \sum_{j=1}^i x_j, 0\le i \le n;\q X_{n+m} := X_n + x + B_{\ch_1-t} + B_{T, T+ \ch_m-\ch_1}, m\ge 1,&\nonumber
\eea
One can easily show that the function $\ol\theta^\e_n$ defined at \eqref{olthen0} also satisfies
\bea
\label{olthen}
&\dis \ol\th^\e_n(\pi_n; t,x) =\sup_{b\in \dbL^0_L(\dbF)} \ol \cE^{L}\Big[e^{\int_t^T\ol b_rdr} \xi(B^\e_{\cdot}(\pi_n,t,x)) +C_0\int_t^Te^{\int_s^T\ol b_rdr} ds\Big],&\\
&\dis \mbox{where}\q \ol b_r := b_{r-t} \1_{[t, \ch_1)}(r) + b_{T+r- \ch_1}   \1_{[\ch_1, T]}(r).&\nonumber
\eea

Fix $(\pi_n, t, x)\in D^\e_{n+1}$. For arbitrary  $(\pi_n',  t',x')  \in D^\e_{n+1}$, define $\ch'_m$,   $\ol b'$, and $N'$ in obvious way.
The advantage of \reff{existence-chn}-\eqref{olthen} is the fact that under this representation one may easily check that, 
\bea
\label{olchm=}
\left.\ba{c}
 \ol b_{t+r} = \ol b_{t'+r},~ 0\le r\le (\ch_1-t) \wedge (\ch_1'-t');\q  \ol b_{\ch_1+r} = \ol b_{\ch_1'+r},~ 0\le r\le (T-\ch_1) \wedge (T-\ch_1');\\
B^\e_{\ch_1, \ch_1+r}(\pi_n,t,x) = B_{\ch_1', \ch_1'+r}^\e(\pi'_n,t',x'),~ r\le (\ch_{N-1}-\ch_1)\wedge (\ch'_{N-1}-\ch_1'),\\
 \ch_{m+1} -  \ch_m = \ch'_{m+1} - \ch_m',\q 1 \le m \le N \wedge N'-2;\\
  N' \le N~\mbox{on}~\{\ch_1 \le \ch_1'\},\q N\le N'~\mbox{on}~\{\ch'_1 \le \ch_1\}.
\ea\right.
\eea
Denote  
\bea
\label{delta1}
\d := \max_{1\le i\le n} [|t_i - t_i'|+|x_i-x_i'|] ~\bigvee ~ [|t-t'|+|x-x'|],
\eea
and for notational simplicity,
\bea
\label{tildeT}
\tilde \ch_1:= \ch_1-t,\q \tilde \ch_1':= \ch_1'-t',\q \tilde T:= T-\ch_1,\q \tilde T':= T-\ch'_1.
\eea
Note that $\ch_1 = (t + \ch^{0,x, \e-L_1(t-t_n)})\wedge T$, $\ch'_1 := (t' + \ch^{0,x', \e-L_1(t'-t'_n)})\wedge T$. By Lemma \ref{lem-cHLip},
\bea
\label{chdifference}
\ol\cE^L[|\ch_1-\ch_1'|] \le C\d,\q\mbox{which implies}\q \sup_{\dbP\in \cP_L} \dbP\big(|\ch_1 - \ch_1'|> \sqrt{\d}\big) \le C\sqrt{\d}.
\eea

Now let $\rho$ denote the modulus of continuity function of $\xi$.
By \reff{olthen} we have
\bea
\label{existence-est1}
\Big| \ol\th^\e_n(\pi_n; t,x) -  \ol\th^\e_n(\pi'_n; t',x')\Big| \le C\sup_{b\in \dbL^0_L(\dbF)} [I_1(b)+I_2(b)+I_3], 
\eea
where
\bea
\label{I}
I_1(b)&:=& \ol\cE^L\Big[\big|e^{\int_t^T \ol b_r dr} -e^{ \int_{t'}^T \ol b'_r dr}\big|\Big];\nonumber\\
I_2(b) &:=& \ol\cE^L\Big[\big|\int_t^Te^{\int_s^T\ol b_rdr} ds - \int_{t'}^Te^{\int_s^T\ol b'_rdr} ds\big|\Big];\\
I_3 &:=& \ol\cE^L\Big[\rho\big(\| B_\cdot^\e(\pi_n,t,x)-B'^\e_\cdot(\pi'_n,t',x') \|_T\big)\Big].\nonumber
\eea
Recall \reff{delta1}. Note that, by \reff{olchm=} and \reff{chdifference},
\beaa
I_1(b) &\le&   C \ol\cE^L\Big[\big|\int_t^T \ol b_r dr - \int_{t'}^T \ol b'_r dr\big|\Big]\\
&=&   C\ol\cE^L\Big[\big|  \int_0^{\tilde \ch_1} \ol b_{t+r} dr -   \int_0^{\tilde \ch'_1} \ol b'_{t'+r} dr + \int_0^{\tilde T} \ol b_{\ch_1+r} dr -\int_0^{\tilde T'} \ol b'_{\ch_1'+r} dr \big|\Big]\\
 &\le& C\ol\cE^L\Big[|\tilde \ch_1 -\tilde \ch'_1| + |\tilde T - \tilde T'| \Big]\le  C\ol\cE^L\Big[|t -t'| + |\ch_1 - \ch_1'| \Big] \le C\d;\\
I_2(b)  &=&   \ol\cE^L\Big[\big|\int_0^{ \tilde \ch_1} e^{\int_{t+s}^T\ol b_rdr} ds - \int_0^{\tilde \ch_1'} e^{\int_{t'+s}^T\ol b'_rdr} ds +  \int_0^{\tilde T} e^{\int_{\ch_1+s}^T\ol b_rdr} ds  - \int_0^{\tilde T'} e^{\int_{ \ch'_1+s}^T\ol b'_rdr} ds\big|\Big]\\
&\le&  C\ol\cE^L\Big[|\int_0^{\tilde \ch_1\wedge \tilde \ch_1'} [e^{\int_{t+s}^T\ol b_rdr}-e^{\int_{t'+s}^T\ol b'_rdr} ]ds + |\tilde \ch_1-\tilde \ch'_1| \\
&&+ \int_0^{\tilde T\wedge \tilde T'} [e^{\int_{\ch_1+s}^T\ol b_rdr} - e^{\int_{ \ch'_1+s}^T\ol b'_rdr}] ds\big| + |\tilde T- \tilde T'|   \Big]\\
&\le& C\d + C\ol\cE^L\Big[\int_0^{\tilde \ch_1\wedge \tilde \ch_1'} |\int_{t+s}^T\ol b_rdr-\int_{t'+s}^T\ol b'_rdr|ds + \int_0^{\tilde T\wedge \tilde T'}|\int_{ \ch_1+s}^T\ol b_rdr - \int_{\ch'_1+s}^T\ol b'_rdr|ds  \Big]\\
&\le&C\d,
\eeaa
where the last estimate follows similar arguments as for $I_1(b)$. Then
\bea
\label{existence-est2}
\Big| \ol\th^\e_n(\pi_n; t,x) -  \ol\th^\e_n(\pi'_n; t',x')\Big| \le C\d + CI_3. 
\eea

The estimate for $I_3$ is more involved.  We first consider the case that $t> t_n$. Denote
\bea
\label{d0}
\d_0 := \min_{1\le i\le n} [t_i - t_{i-1}] \wedge [t-t_n] > 0,
\eea
and let $\d_1 \le {1\over 2}\d_0$ which will be specified later. Consider $\d \le \d_1$, then 
\beaa
 \min_{1\le i\le n} [t'_i - t'_{i-1}] \wedge [t'-t'_n] \ge \d_1.
\eeaa
 For any  $m\ge 1$,  by Lemma \ref{lem-chn} we see that
\bea
\label{existence-chN}
&\dis\sup_{\dbP\in \cP_L} \dbP(\O_m^c) \le {C_\e\over m},&\\
&\dis\mbox{where}~\O_m :=  \{N \le m\} \cap \{N' \le m\}\cap \big\{\sup_{0\le s_1 < s_2\le 2T} {|B_{s_1,s_2}|\over (s_2-s_1)^{1\over 3}} \le m\big\}\cap \{\|B\|_{2T} \le m\},\nonumber
\eea
and the constant $C_\e$ is independent of $(\pi_n, t,x)$ and $(\pi_n', t',x')$. Moreover, by \reff{chlowerbound} we have
\bea
\label{chdiff}
\sup_{\dbP\in \cP^L}\Big[ \dbP\big(\ch_{i+1} - \ch_i < \d_1,  \ch_{i+1} < T
\big) + \dbP\big(\ch'_{i+1} -\ch'_i < \d_1, \ch'_{i+1}<T
\big)\Big] \le C_\e\d_1,\q i\ge 1.
\eea
This implies 
\beaa
\sup_{\dbP\in \cP^L}\dbP\Big(\big[\bigcup_{i=1}^{N-2}\{\ch_{i+1} - \ch_i < \d_1\}\cap \{N \le m\}\big] \bigcup \big[\bigcup_{i=1}^{N'-2}\{\ch'_{i+1} - \ch'_i < \d_1\}\cap \{N' \le m\}\big]\Big) \le C_\e m\d_1.
\eeaa
Thus, by \reff{chdifference} and assuming $\sqrt{\d}\le \d_1$,
\bea
\label{existence-d1}
&\dis \sup_{\dbP\in \cP^L} \dbP(\O_{m,\d_1}^c) \le C_\e\big[{1\over m} +   m \d_1], \q\mbox{where}&\\
 &\dis\O_{m,\d_1} := \O_m \bigcap \{|\ch_1-\ch_1'|\le \sqrt{\d}\} \bigcap \Big(\bigcap_{i=1}^{N-2}\{\ch_{i+1} - \ch_i \ge \d_1\}\Big) \bigcap \Big(\bigcap_{i=1}^{N'-2}\{\ch'_{i+1} - \ch'_i \ge \d_1\}\Big).&\nonumber
\eea
Denote $\t_0:=0$,  $\t_i:= t_i$,   $1\le i\le n$, $\t_{n+i} := \ch_i$, $1\le i\le m$, and define $\t'_i$ similarly.  Note that $t-t_n \ge \d_1$, $t'-t_n'\ge \d_1$, and $\ch_i- \ch_i' = \ch_1-\ch_1'$, $1\le i< N\wedge N'$, thanks to the third line of \reff{olchm=}.  Then, assuming  $m \ge \sum_{i=1}^n[ |x_i| +|x_i'|] +|x|+|x'| +n$,  on $\O_{m,\d_1}$ it holds that
\bea
\label{Omd1}
\left.\ba{c}
\dis \inf_{1\le i< n+N} [\t_{i} - \t_{i-1}] \ge \d_1, \q \inf_{1\le i< n+N'} [\t'_{i} - \t'_{i-1}] \ge \d_1,\q  \sup_{1\le i< n+N\wedge N'}  |\t_i-\t_i'|\le \sqrt{\d};\\
\dis\sup_{1\le i\le N}|X_i|  \le  2m,\q  \sup_{1\le i\le N'}|X_i'|  \le  2m,\q \sup_{1\le i < n+N\wedge N'} |X_i-X_i'| \le Cm \d^{1\over 6}.
\ea\right.
\eea

We now estimate $I^3_s:=|B_s^\e(\pi_n, t, x)- B_s^\e(\pi_n',t',x')|$ on the set $\O_{m,\d_1}$. Without loss of generality we assume  $\sqrt{\d} \le {\d_1\over 2}$, $\ch_1' \le \ch_1$ and thus  $N \le N'$. We estimate $I^3_s$ in several cases.

{\it Case 1.} $s \in [\t_i + \sqrt{\d}, \t_{i+1}-\sqrt{\d}]$ for some $0\le i\le n+N-2$. Then, noting that $|\t_j-\t'_j|\le \sqrt{\d}$ for $j=i,i+1$,  we have $s\in [\t_i', \t_{i+1}']$ and thus
\bea
\label{Case1}
I^3_s &=& \Big|X_i + {s-\t_i\over \t_{i+1}-\t_i}[X_{i+1}-X_i] - X'_i - {s-\t'_i\over \t'_{i+1}-\t'_i}[X'_{i+1}-X'_i]\Big|\nonumber\\
&=& \Big|[X_i-X_i'] + {s-\t'_i\over \t'_{i+1}-\t'_i}[(X_{i+1}-X_i)-(X_{i+1}'-X_i')] \\
&&+\big[[{s-\t_i\over \t_{i+1}-\t_i} - {s-\t_i'\over \t_{i+1}-\t_i}] +[{s-\t_i'\over \t_{i+1}-\t_i}- {s-\t'_i\over \t'_{i+1}-\t'_i}]\big][X_{i+1}-X_i] \Big|\nonumber.
\eea
Now by \reff{Omd1} we can easily see that $I^3_s\le Cm\d^{1\over 6} + {Cm\over \d_1}\sqrt{\d}$.

{\it Case 2.} $|s-\t_i|\le \sqrt{\d}$ for some $0\le i <n+N$. Then $|s-\t_i'|\le 2\sqrt{\d}$, and
\bea
\label{Case2}
I^3_s &\le& |B^\e_s(\pi_n,t,x) - B^\e_{\t_i}(\pi_n,t,x)| + |X_i-X_i'|+|B^\e_{\t_i'}(\pi'_n,t',x') - B^\e_s(\pi_n',t',x')|.
\eea

{\it Case 2.1.} Assume $0\le i\le n$. Then  \reff{Case2} and \reff{Omd1} lead to 
\beaa
I^3_s&\le& {Cm\over \d_1} [|s-\t_i|+|s-\t'_i|] + Cm\d^{1\over 6} \le Cm\d^{1\over 6} + {Cm\over \d_1}\sqrt{\d}
\eeaa

{\it Case 2.2.} Assume $n< i < n+N$.  Note that, when $s > \t_i$,
\beaa
|B^\e_{\t_i,s}(\pi_n,t,x) | = {s-\t_i\over \t_{i+1}-\t_i} |X_{i+1}-X_{i}|\le {s-\t_i\over \t_{i+1}-\t_i}  m  (\t_{i+1}-\t_i)^{1\over 3} \le m (s-\t_i)^{1\over 3} \le m \d^{1\over 6}.
\eeaa
Similarly we have the other related estimates. Then \reff{Case2} and \reff{Omd1} lead to $I^3_s\le Cm\d^{1\over 6}$.

{\it Case 3.} $s\in [\t_{n+N-1} + \sqrt{\d}, T]$ and $T-\t_{n+N-1}\le\d^{1\over 4}$. Then $s-\t_{n+N-1} \le \d^{1\over 4}$ and $0\le s - \t'_{n+N-1} \le \d^{1\over 4} + \sqrt{\d}\le 2 \d^{1\over 4}$. When $N'=N$, then of course $s \in [\t_{n+N-1}', \t_{n+N}']$. When $N<N'$, then $\t'_{n+N} - \t'_{n+N'-1} \ge \d_1$ and thus we still have  $s \in [\t_{n+N-1}', \t_{n+N}']$ whenever $2 \d^{1\over 4}\le \d_1$. Now following the arguments in Case 2.2, by \reff{Case2}  and \reff{Omd1} we can easily see that $I^3_s \le   Cm\d^{1\over 12}$.

{\it Case 4.} $s\in [\t_{n+N-1} + \sqrt{\d}, T]$ and $T-\t_{n+N-1}\ge\d^{1\over 4}$. Denote  $i:= n+N-1$ for notational simplicity, then $\t_{i+1}=T$.  Similar to \reff{olchm=}, by \reff{existence-chn} one can easily see that $\t'_{i+1} - \t'_i\ge \t_{i+1}-\t_i$. Then, together with \reff{Omd1} we have
\bea
\label{tN}
\t'_{i+1} - \t'_{i} \ge  T-\t_{i} \ge \d^{1\over 4},\q \t_{i+1}-\t'_{i+1} \le  \t_{i} -\t'_{i} \le \sqrt{\d},\q
 \t'_{i+1}-\t'_{i} \le  T-\t_{i} + \sqrt{\d}.
\eea

{\it Case 4.1.} Assume $s \in [\t_i', \t_{i+1}']$.  Then, by \reff{Case1}, \reff{Omd1}, and the first two inequalities of \reff{tN}, we have 
$I^3_s \le Cm\d^{1\over 6} + {Cm\over \d^{1\over 4}}\sqrt{\d} \le Cm \d^{1\over 4}$.

{\it Case 4.2.} Assume $s \notin [\t_i', \t_{i+1}']$.  In this case, we must have $N'>N$ and thus $\t'_{i+1} - \t'_i\ge \d_1$.  By the second inequality of \reff{tN} we see that $T-\t'_{i+1} \le \sqrt{\d} < \d_1$, then  we must have $N'= N+1$ and thus $s\in [\t'_{i+1}, \t'_{i+2}]$. By \reff{Omd1} and  \reff{tN} again we can see that 
\beaa
 T- \t_{i} \ge \d_1 -\sqrt{\d},\q 0\le s - \t'_{i+1} \le T-\t'_{i+1}\le \sqrt{\d},\q T-s \le T-\t'_{i+1} \le \sqrt{\d},&
\eeaa
Then, recalling the notations in \reff{tildeT},
\beaa
I^3_s &=& \Big|X_{i} + B_{\tilde T+\t_i, \tilde T+T} - {T-s\over T-\t_i} B_{\tilde T+\t_i, \tilde T+T} - X'_{i}-B_{\tilde T'+\t'_i, \tilde T'+\t'_{i+1}}  - {s-\t'_{i+1}\over T-\t'_{i+1}} B_{\tilde T'+\t'_{i+1}, \tilde T' + T} \Big|\\
&\le& Cm\d^{1\over 6} + {T-s\over T-\t_i} |B_{\tilde T+\t_i, \tilde T+T}| + |B_{\tilde T+\t_i, \tilde T+T} -B_{\tilde T'+\t'_i, \tilde T'+\t'_{i+1}}| + Cm |T-\t'_{i+1}|^{1\over 3}\\
&\le& Cm\d^{1\over 6} + Cm {\sqrt{\d}\over \d_1 -\sqrt{\d}} + C[|\tilde T- \tilde T'|^{1\over 3} + |\t_i - \t_i'|^{1\over 3} + |T-\t'_{i+1}|^{1\over 3}] \le    Cm\d^{1\over 6} + Cm {\sqrt{\d}\over \d_1}. 
\eeaa

Put all the cases together, we have
\beaa
I^3_s \le  Cm\d^{1\over 12} + Cm {\sqrt{\d}\over \d_1} \q\mbox{on} ~\O_{m,\d_1} \cap \{\ch_1' \le \ch_1\}.
\eeaa
We may get the same estimate on $\O_{m,\d_1} \cap \{\ch_1\le \ch_1'\}$. Plug this and \reff{existence-d1} into \reff{existence-est2}, we obtain
\beaa
\Big| \ol\th^\e_n(\pi_n; t,x) -  \ol\th^\e_n(\pi'_n; t',x')\Big| \le C\d + {C_\e\over m} +C_\e m\d_1+ Cm\d^{1\over 12} + Cm {\sqrt{\d}\over \d_1}.
\eeaa
Set $\d_1 := {1\over m^2}\wedge {\d_0\over 2}$. Then, whenever $\d \le {1\over m^{48}}$, 
\beaa
\Big| \ol\th^\e_n(\pi_n; t,x) -  \ol\th^\e_n(\pi'_n; t',x')\Big| \le  {C_\e\over m} + C_\e m^3\d^{1\over 12} \le {C_\e\over m}.
\eeaa
Since $m\ge 1$ is arbitrary, we see that $\ol\th^\e_n$ is continuous at $(\pi_n, t, x)$, in the case that $t> t_n$.

In the case $t=t_n$ and thus $x=0$, we modify \reff{d0} as
\beaa
\d_0 := \min_{1\le i\le n} [t_i - t_{i-1}] > 0.
\eeaa
Note that, in this case $|x'|\le \d$ and one can easily see that \reff{chdiff} still holds true for $i=0$. Then following almost the same arguments as in previous case we may prove that $\ol\th^\e_n$ is continuous at $(\pi_n, t_n, 0)$.
\qed

\end{document}